\documentclass[leqno, 12pt]{article}

\usepackage{amsmath}
\usepackage{amsfonts}
\usepackage{amssymb}
\usepackage{amsthm} 
\usepackage{graphicx}

\usepackage{scalerel,stackengine}
\stackMath
\newcommand\widecheck[1]{%
\savestack{\tmpbox}{\stretchto{%
  \scaleto{%
    \scalerel*[\widthof{\ensuremath{#1}}]{\kern-.6pt\bigwedge\kern-.6pt}%
    {\rule[-\textheight/2]{1ex}{\textheight}}%WIDTH-LIMITED BIG WEDGE
  }{\textheight}% 
}{0.5ex}}%
\stackon[1pt]{#1}{\scalebox{-1}{\tmpbox}}%
}
\usepackage{scalerel,stackengine}
\stackMath
\newcommand\widebreve[1]{%
\savestack{\tmpbox}{\stretchto{%
  \scaleto{%
    \scalerel*[\widthof{\ensuremath{#1}}]{\kern-.1pt\bigcup\kern-.1pt}  
    {\rule[-\textheight/2]{1ex}{\textheight}}%WIDTH-LIMITED BIG CUP
  }{\textheight}% 
}{0.5ex}}%
\stackon[1pt]{#1}{\scalebox{1}{\tmpbox}}%
}
\usepackage{comment}

\def\toL1{\buildrel \mathit{L_1}\over\longrightarrow}

\def\towpr1{\buildrel w.pr.1\over\to}
\def\11{\buildrel 1-1\over\longleftrightarrow}

\def\~{\tilde}
\def\stms{\!\setminus\!} 
\def\^{\hat}
\def\u{\breve}
\def\c{\check}

\def\implies{\Rightarrow}

\def\Sig{\Sigma}
\def\sig{\sigma}
\def\Lam{\Lambda}
\def\Gam{\Gamma}
\def\gam{\gamma}
\def\lam{\lambda}
\def\alp{\alpha}

\def\del{\delta}

\def\eps{\epsilon}

\def\ome{\omega}
\def\Ome{\Omega}

\def\E{{\rm E}}

\def\R1{{\bf R}^1}
\def\B1{{\cal B}^1}

\def\nid{\noindent}

\def\f12{\frac{1}{2}}

\def\ts{\textstyle}
\def\t12{\textstyle{1\over2}}

\def\smalltype{\let\rm=\eightrm \let\bf=\eightbf \let\it=\eightit \let\sl=\eightsl
 \baselineskip=9pt \rm}
\font\tenrm=cmr10
\font\tenbf=cmbx10
\font\tenit=cmti10
\font\tensl=cmsl10
\def\medtype{\let\rm=\tenrm \let\bf=\tenbf \let\it=\tenit \let\sl=\tensl
 \baselineskip=12pt \rm}
\font\twelverm=cmr12
\font\twelvebf=cmbx12
\font\twelveit=cmti12 
\font\twelvesl=cmsl12
\def\bigtype{\let\rm=\twelverm \let\bf=\twelvebf \let\it=\twelveit \let\sl=\twelvesl
 \baselineskip=14pt \rm}
 
\begin{document}

%%COMMENT
%\begin{comment}
\title{On the Feasibility of Parsimonious Variable Selection for Hotelling's $T^2$-test}
\author{Michael D. Perlman\footnote{mdperlma@uw.edu.}\\
University of Washington}

\maketitle

\begin{abstract}

\nid  Hotelling's $T^2$-test for the mean of a multivariate normal distribution is one of the triumphs of classical multivariate analysis. It is uniformly most powerful among invariant tests, and admissible, proper Bayes, and locally and asymptotically minimax among all tests. Nonetheless, investigators often prefer non-invariant tests, especially those obtained by selecting only a small subset of variables from which the $T^2$-statistic is to be calculated, because such reduced statistics are more easily interpretable for their specific application. Thus it is relevant to ask the extent to which power is lost when variable selection is limited to very small subsets of variables, e.g. of size one (yielding univariate Student-$t^2$ tests) or size two (yielding bivariate $T^2$-tests).  This study presents some evidence, admittedly fragmentary and incomplete, suggesting that in some cases no power may be lost over a wide range of alternatives.

%linear combinations of small subsets of the individual variables preferred to general linear combinations.....

\end{abstract}
%\end{comment}
%%END COMMENT

\vskip4pt

\newpage

{\it This work is dedicated to the memory of Somesh Das Gupta, my colleague, teacher, and friend.}
\vskip4pt

\nid{\bf 1. Introduction.} This study is motivated by a re-examination of the variable-selection problem for Hotelling's $T^2$-test (closely related to variable selection for linear discriminant analysis). After some notational preliminaries in \S1.1, Hotelling's $T^2$ is reviewed in \S1.2. The variable-selection problem is described in \S1.3, where the substance of this investigation is described. 

\nid{\it 1.1. The noncentral $f$-distribution.} Let $\chi_m^2(\lam)$ denote a noncentral chi-square random variable with $m$ degrees of freedom and noncentrality parameter $\lam>0$. The noncentral $f_{m,n}(\lam)$ distribution (nonnormalized) with $m$ and $n$ degrees of freedom and noncentrality parameter $\lam>0$ is the distribution of the ratio $\chi_m^2(\lam)/\chi_n^2$ (also denoted by $f_{m,n}(\lam)$), where the numerator and denominator are independent chi-square random variables and $\chi_n^2\equiv\chi_n^2(0)$. The upper $\alp$-quantile of $f_{m,n}\equiv f_{m,n}(0)$ is denoted by $f_{m,n}^\alp$, so that 
\begin{equation}\label{fquantile}
\Pr[f_{m,n}>f_{m,n}^\alp]=\alp.
\end{equation}

The {\it noncentral $f_{m,n}$-test of size} $\alp\ge0$ for the problem of testing $\lam=0$ vs. $\lam>0$ has power function given by
\begin{align}
\pi_\alp(\lam;m,n)&=\Pr[f_{m,n}(\lam)>f_{m,n}^\alp]\label{pialp}\\
 &=e^{-\frac{\lam}{2}}\sum_{k=0}^\infty\ts(\frac{\lam}{2})^k\frac{1}{k!}c_{m,n;k;\alp},\label{6}\\
 c_{m,n;k;\alp}&\equiv\Pr[f_{m+2k,n}>f_{m,n}^\alp];\label{cmnk}
\end{align}
see Das Gupta and Perlman [DGP] (1974), eqn.(2.1). Clearly $\pi_\alp(\lam;m,n)$ is decreasing in $\alp$, with $\pi_0(\lam;m,n)=0$. Because  $f_{m,n}(\lam)$ has strictly monotone likelihood ratio in $\lam$, $\pi_\alp(\lam;m,n)$ is strictly increasing in $\lam$.

It will be convenient to work with the (central) {\it beta distribution} $b_{m,n}$:
\begin{equation}\label{beta}
b_{m,n}:=\frac{f_{m,n}}{f_{m,n}+1}\equiv\frac{\chi_m^2}{\chi_m^2+\chi_n^2};
\end{equation}
clearly, $b_{m,n}= 1-b_{n,m}$. The upper and lower $\alp$-quantiles of $b_{m,n}$ are denoted by $b_{m,n}^\alp$ and $b_{m,n;\alp}$, respectively, so that
\begin{equation}\label{betafact}
b_{m,n}^\alp= 1-b_{n,m;\alp}.
\end{equation}
Thus by \eqref{cmnk},
\begin{equation}\label{cmnk2}
 c_{m,n;k;\alp}=\Pr[b_{n,m+2k}<b_{n,m;\alp}].
 \end{equation}
 
The probability density function of $b_{m,n}$ is given by
\begin{equation}\label{betapdf}
\phi_{m,n}(b)\equiv\frac{\Gam(\frac{m+n}{2})}{\Gam(\frac{m}{2})\Gam(\frac{n}{2})}b^{\frac{m}{2}-1}(1-b)^{\frac{n}{2}-1},\qquad 0<b<1.
\end{equation}
%\vskip1pt

\nid{\it 1.2. Hotelling's $T^2${-}test.} Let $X_i:p\times 1$, $i=1,\dots,N$ ($N\ge p+1$) be a random sample  from the $p$-dimensional multivariate normal distribution $N_p(\mu,\Sig)$, where $\mu\ (p\times1)\equiv(\mu_1,\dots,\mu_p)' \in\mathbb{R}^p$ and $\Sig\ (p\times p)\equiv(\sig_{ij})$ is positive definite. The problem of testing 
\begin{equation}\label{H0K}
H_0:\mu=0\quad\text{vs.}\quad K:\mu\ne0
\end{equation}
with $\Sig$ unknown is invariant under the group action $X_i\to A X_i$, $i=1,\dots,N$, where $A\in GL(p)$, the group of all nonsingular $p\times p$ matrices. The maximal invariant statistic under $GL(p)$ is given by Hotelling's $T^2$ statistic:
\begin{equation}\label{T2}
T^2:= N\bar X'S^{-1}\bar X,
\end{equation}
where $\bar X=\sum_{i=1}^NX_i$ and $S=\sum_{i=1}^N(X_i-\bar X)(X_i-\bar X)'$. Its distribution is 
\begin{equation}\label{T2A}
T^2\sim f_{p,N-p}(\Lam),
\end{equation}
where
\begin{equation}\label{lammuSig}
\Lam=N\mu'\Sig^{-1}\mu,
\end{equation}
the maximal invariant parameter. Therefore the uniformly most powerful invariant size-$\alp$ test rejects $H_0$ if $T^2>f_{p,N-p;\alp}$, with power function $\pi_\alp(\Lam;p,N-p)$; cf.  [A] Theorem 5.6.1).\footnote{\label{FootAdmiss} Among the larger class of all tests, invariant or non-invariant, the $T^2$ test is admissible [S], proper Bayes [KS], and locally and asymptotically minimax for small and large values of $\Lam$, respectively [GK].}

It is useful to express $\Lam$ in terms of scale-free parameters, that is,
\begin{equation}\label{lamgamR}
\Lam=N\gam'R^{-1}\gam,
\end{equation}
where $R\equiv(\rho_{ij})$ is the $p\times p$ correlation matrix determined by $\Sig$ and
\begin{align}
\ts\gam\ (p\times1)\equiv(\gam_1,\dots,\gam_p)':=\big(\frac{\mu_1}{\sqrt{\sig_{11}}},\dots,\frac{\mu_p}{\sqrt{\sig_{pp}}}\big)'.\label{gamma}
\end{align}
The testing problem \eqref{H0K} can be stated equivalently as that of testing
\begin{equation}\label{H0KR}
H_0:\gam=0\quad\text{vs.}\quad K:\gam\ne0
\end{equation}
with $R$ unknown.
\vskip4pt

\nid{\it 1.3. The $T^2$ variable-selection problem.}  Denote the components of $\bar X$ by $\bar X_j$, $j=1,\dots,p$, and those of $S$ by $s_{jk}$, $j,k=1,\dots,p$. Let $\Ome_p$ be the collection of all nonempty subsets of $\{1,\dots,p\}$. For $\ome\in\Ome_p$ denote the $\ome$-subvector of $\bar X$ by $\bar X_\ome$, the $\ome$-submatrix of $S$ by $S_\ome$, and similarly define $\gam_\ome$ and $R_\ome$. The $T^2$-statistic based on $(\bar X_\ome,S_\ome)$ is given by
\begin{align}
T_\ome^2&\equiv N\bar X_\ome'S_\ome^{-1}\bar X_\ome\sim f_{|\ome|,\,N-|\ome|}(\Lam_\ome),\label{T2ome}\\
\Lam_\ome&\equiv\Lam_\ome(\gam,R)\ =\ N\gam_\ome'R_\ome^{-1}\gam_\ome.\label{lamome}
\end{align}
($T_{\Ome_p}^2=T^2$, $\Lam_{\Ome_p}=\Lam\equiv\Lam(\gam,R)$.) The test that rejects $H_0$ if $T_\ome^2>f_{|\ome|,N-|\ome|;\alp}$ has size $\alp$ for $H_0$, and its power function is given by 
\begin{equation}\label{pilamome}
\pi_\alp(\Lam_\ome;\,|\ome|,\,N-|\ome|).
\end{equation}
This $T_\ome^2$-test is not invariant under $GL(p)$ but it is admissible for testing $H_0$ vs. $K$, being a unique proper Bayes test for a prior distribution that, under $K$, assigns mass 1 to $\{\mu\mid \mu_{\Ome_p\setminus\ome}=0\}$; cf. [KS], [MP]. 

The $T^2$ variable-selection problem is that of choosing a parsimonious subset $\ome$ such that the $T_\ome^2$-test maintains high power against alternatives deemed likely in some sense. Because $(\gam,R)$ is unknown, variable selection in practice is traditionally approached by forward and/or backward selection procedures based on a preliminary sample that yields estimates of $(\gam,R)$; see the Appendix. In general, all $2^p-1$ nonempty subsets $\ome$ must be considered. 

Recently I consulted on such a variable-selection problem. The investigator, a web-page designer, had observed 20 physiological variables (blood pressure, temperature, heart rate, etc.) on each of 100 subjects (the numbers are approximate). He wished to compare their responses to a new web-page design with their responses to the current design. The overall $T^2$-statistic, based on a linear combination of all 20 variables, indicated a significant difference between the two sets of responses. However, the client desired to find a more readily interpretable measure of difference, namely a $T_\ome^2$-statistic based on a very small subset $\ome$ of the 20 variables, hopefully with $|\ome|=1$ or 2. Such a desire is not atypical of investigators presented with a multivariate data analysis. Therefore it occurred to me to wonder just how much power would be lost by restricting variable selection to small subsets $\ome$, for example to single variables or pairs of variables. 

To state this more precisely, define
\begin{equation}\label{omestar}
\^\ome_\alp(\gam,R)=\underset{\ome\in\Ome_p}{\arg\max}\ \pi_\alp(\Lam_\ome(\gam,R);\,|\ome|,\,N-|\ome|).
\end{equation}
Thus $\^\ome_\alp(\gam,R)$ is the (not necessarily unique) subset $\ome$ of variables that maximizes the power of the size-$\alp$ $T_\ome^2$-test to detect the alternative $(\gam,R)$ if the actual value of $(\gam,R)$ were revealed by an oracle. Whereas the admissibility of the overall size-$\alp$ $T^2$-test dictates that its power cannot be dominated by that of the size-$\alp$ $T_\ome^2$-test when $\ome\ne\Ome_p$,
% for any proper subset $\ome\subset\Ome_p$, 
might it happen that $|\^\ome_\alp(\gam,R)|$ is small, perhaps 1 or 2, over a fairly wide range of parameter values $(\gam,R)$?  If so, then one might with some confidence limit variable selection to consideration of single variables (univariate $t^2$-tests) or pairs of variables (bivariate $T^2$-tests) -- thus the overall $p$-variate $T^2$-test would be ruled out {\it a priori}.

Of course, such a radical suggestion flies in the face of 100 years of multivariate statistical theory, of which I have been but one of many proponents. Nonetheless, this report presents some evidence, admittedly fragentary and incomplete, indicating that this suggestion might not be entirely inappropriate in applications. In Sections 2, 3, and 4, a few very special cases are considered where tractable algebraic expressions for the asymptotic ($\Lam_\ome\to\infty$), local ($\Lam_\ome\to0$), and/or exact values of $\pi_\alp(\Lam_\ome;|\ome|, N-|\ome|)$ are available. These in turn can be utilized to compare the powers of $T_\ome^2$ and $T^2$. 

Examples 2.1 and 3.1 treat only the simplest possible case: the bivariate case ($p=2$) with $N=3$.\footnote{\label{Foot4} However, Giri, Kiefer, and Stein [GKS] also began their study of Hotelling's $T^2$ test by considering only this simplest case $p=2$, $N=3$.} Here it is shown that $|\^\ome_\alp(\gam,R)|=1$ over large portions of the asymptotic and local regions of the alternative hypothesis $K$. This implies that the power of at least one of the two univariate Student $t^2$-tests ($|\ome|=1$) exceeds that of the overall (bivariate) $T^2$-test for most alternatives $(\gam,R)$ in these regions.

In Example 4.4 this result is extended to the entire alternative hypothesis $K$, both for $N=3$ and $N=5$, but only under the highly restrictive and impracticably vague condition that $\alp$ be sufficiently small, with ``sufficiently small" determined by the value of the unknown noncentrality parameter -- see \eqref{ineq15} and \eqref{ineq155}. 

Examples 2.2 and 3.2 go beyond the bivariate case. Here $p\ge3$, $N=p+2$, and the powers of all possible bivariate $T_\ome^2$-tests ($|\ome|=2$) are compared to the power of the overall $p$-variate $T^2$ test, again only for asymptotic and local alternatives and, furthermore, only for very special configurations of $\gam$ and $R$. In these cases, admittedly highly restrictive, the bivariate  $T_\ome^2$-tests dominate the overall $T^2$ over a substantial portion of the alternative hypothesis $K$. This does not establish that $|\^\ome_\alp(\gam,R)|=2$ but again suggests that variable search might be limited to small variable subsets $\ome$. 

Of course, a much more comprehensive study is needed to confirm the efficacy of such an approach to variable selection. It is hoped that this report might encourage such an investigation.
\vskip6pt

%\newpage

\nid{\bf 2. Some asymptotic power comparisons.} The power function of the $T_\ome^2$-test is 
\begin{equation}\label{pinew}
\pi_\alp(\Lam_\ome)\equiv\pi_\alp(\Lam_\ome;|\ome|,N-|\ome|)
\end{equation}
(recall \eqref{T2ome}-\eqref{lamome}). It follows from eqn. (3.4) in [MP] that as $\Lam_\ome\to\infty$, 
\begin{align}
\pi_\alp(\Lam_\ome)&\sim \ts1-\exp[{-\frac{\Lam_\ome}{2}(f_{|\ome|,N-|\ome|}^\alp+1})^{-1}]\nonumber\\
&=\ts1-\exp[-\frac{\Lam_\ome}{2}b_{N-|\ome|,|\ome|;\alp}].\label{asymp1}
%&=1-e^{-\frac{\Lam_\ome}{2}}b_{N-|\ome|,|\ome|;\alp}.%\quad\mathrm{as}\ \Lam_\ome\to\infty.%\mathrm{ess\ inf}\{\}
\end{align}
Thus for two subsets $\ome,\,\ome'$ s.t. $\ome\subset\ome'$,
%$|\ome|\subset|\ome'|$,  
$\exists\;\Lam_{|\ome|,|\ome'|,N;\alp}^*>0$ s.t.
%$\exists\;\Lam_\alp\equiv\Lam_{|\ome|,|\ome'|;\alp}>0$ s.t.
\begin{equation}\label{asymp3}
\Lam_\ome>\max\left(\Lam_{|\ome|,|\ome'|,N;\alp}^*,\;\ts\frac{b_{N-|\ome'|,|\ome'|;\alp}}{b_{N-|\ome|,|\ome|;\alp}}\Lam_{\ome'}\right)\ \implies\ \pi_\alp(\Lam_\ome)>\pi_\alp(\Lam_{\ome'}).
\end{equation}
Therefore power comparisons of $T_\ome^2$ and $T_{\ome'}^2$ for distant alternatives\footnote{\label{Foot2} We do not claim to know the values of $\Lam_{|\ome|,|\ome'|,N;\alp}^*$, even approximately.} require determination of the lower quantiles $b_{n,m;\alp}$. This can be done explicitly in Examples 2.1 and 2.2 below. Although these examples are of very limited scope\footnote{\label{Foot7} But see Footnote \ref{Foot4}.} they begin to suggest that variable subset selection sometimes can be limited to very small subsets $\ome\in\Ome_p$, e.g., singletons in the bivariate Example 2.1, or pairs (including singletons) in Example 2.2. 
%Thus for two subsets $|\ome|\subset|\ome'|$, as $\lam_\ome\to\infty$,
%\begin{equation}\label{asymp2}
%\pi_\alp(\lam_\ome)>\pi_\alp(\lam_{\ome'})\iff \lam_\ome b_{N-|\ome|,|\ome|;\alp}]>\lam_{\ome'}b_{N-|\ome'|,|\ome'|;\alp}.
%\end{equation}

To simplify the notation, set 
\begin{equation}\label{Qdef}
Q_{|\ome|,|\ome'|,N;\alp}:=\frac{b_{N-|\ome'|,|\ome'|;\alp}}{b_{N-|\ome|,|\ome|;\alp}}<1.
\end{equation}
The quantile $b_{n,m;\alp}$ satisfies
\begin{equation}\label{balpha}
\alp=\ts\frac{\Gam\frac{n+m}{2})}{\Gam(\frac{n}{2})\Gam(\frac{m}{2})}\int_0^{b_{n,m;\alp}}b^{\frac{n}{2}-1}(1-b)^{\frac{m}{2}-1}db
\end{equation}
For the simple cases $n=2$ or $m=2$,
\begin{equation}\label{twobs}
b_{2,m;\alp}=1-(1-\alp)^{\frac{2}{m}},\quad b_{n,2;\alp}=\alp^{\frac{2}{n}}.
\end{equation}

%\newpage

\nid{\bf Example 2.1.} In the bivariate case $p=2$, abbreviate the singleton subsets $\{1\}$ and $\{2\}$ of $\Ome_2$ by $1$ and $2$ respectively. 
We shall compare the powers $\pi_\alp(\Lam_1;\,1,\,2)$ and $\pi_\alp(\Lam_2;\,1,\,2)$ of the two univariate size-$\alp$ $t^2$-tests to the power $\pi_\alp(\Lam;\,2,\,1)$ of  the overall (bivariate) size-$\alp$  $T^2$-test for distant alternatives. 

Assume that $\gam_1\ne0$ (recall \eqref{gamma}) and set  
\begin{align}
%\mu&=\begin{pmatrix}\mu_1\\\mu_2\end{pmatrix},\quad\Sig=\begin{pmatrix}\sig_1^2&\sig_1\sig_2\rho\\ \sig_1\sig_2\rho&\sig_2^2\end{pmatrix},\\
\ts\eta=\frac{\gam_2}{\gam_1},\quad\ \ \,\rho=\rho_{12},\label{etarho}
\end{align}
where $-1<\rho<1$ , so by \eqref{lamgamR} and \eqref{lamome},
%$\mu=(\mu_1,\mu_2)'$ and $\Sig=\begin{pmatrix}\sig_1^2&\sig_1\sig_2\rho\\ \sig_1\sig_2\rho&\sig_2^2\end{pmatrix}$, then
\begin{equation}\label{lams}
\ts\Lam_1=N\gam_1^2,\quad\Lam_2=N\eta^2\gam_1^2,\quad \Lam=N\big(\frac{1-2\eta\rho+\eta^2}{1-\rho^2}\big)\gam_1^2.
\end{equation}
Without loss of generality we can assume that $|\gam_1|\ge|\gam_2|$, so $0\le\eta^2\le1$ and
$\max(\Lam_1,\Lam_2)=\gam_1^2$. Thus the alternative hypotheses $K$ can be represented as 
\begin{equation}\label{altgametarho}
K=\{(\gam_1,\eta,\rho)\mid |\gam_1|>0,\,|\eta|\le1,\,|\rho|<1\},
\end{equation}
while $\^\ome_\alp(\gam,R)$ can be re-expressed as $\^\ome_\alp(\gam_1,\eta,\rho)$.
%\begin{equation}\label{hatomegametarho}
%\^\ome_\alp(\gam_1,\eta,\rho).
%\end{equation}

Because $\max(\Lam_1,\Lam_2)\le\Lam$, it follows from \eqref{asymp3} and \eqref{lams} that
%for $j=1,2$,
\begin{align}
&\ts\qquad\gam_1^2>\max\left(\frac{1}{N}\Lam_{1,2,N;\alp}^*,\;\ts Q_{1,2,N;\alp}\big(\frac{1-2\eta\rho+\eta^2}{1-\rho^2}\big)\gam_1^2\right)\label{asymp4}\\
\implies\ &\max(\pi_\alp(\Lam_1;1,2),\pi_\alp(\Lam_2;1,2))>\pi_\alp(\Lam;2,1).\label{asymp4a}\\
\implies\ & |\^\ome_\alp(\gam_1,\eta,\rho)|=1.\label{hatome1}
\end{align}

In the simplest case $N=3$, \eqref{twobs} yields the explicit expression
\begin{equation}\label{Q12}
\ts Q_{1,2,3;\alp}=\frac{b_{1,2;\alp}}{b_{2,1;\alp}}=\frac{\alp}{2-\alp}, while
%=\frac{\alp^2}{1-(1-\alp)^2},
\end{equation}
the inequality in \eqref{asymp4} is equivalent to
%\begin{equation}\label{asympt5}
%\ts1>\max\left(\frac{\lam_{1,2,3;\alp}^*}{\gam_1^2},\;\ts\bet_\alp\big(\frac{1-2\eta\rho+\eta^2}{1-\rho^2}\big)\right).
%\end{equation}
%???If, without loss of generality, we assume that $|\gam_1|\ge|\gam_2|$, so $0\le\eta^2\le1$, then \eqref{asympt5} becomes
\begin{equation}\label{asympt6}
\ts1>\max\left(\frac{\Lam_{1,2,3;\alp}^*}{3\gam_1^2},\;\ts Q_{1,2,3;\alp}\big(\frac{1-2\eta\rho+\eta^2}{1-\rho^2}\big)\right).
\end{equation}
Note that
\begin{align}
&\ts1>Q_{1,2,3;\alp}\big(\frac{1-2\eta\rho+\eta^2}{1-\rho^2}\big)\\
\iff 0>\rho^2-&2Q_{1,2,3;\alp}\eta\rho+Q_{1,2,3;\alp}(1+\eta^2)-1=:h_{\alp,\eta}(\rho).\label{lamalp}
\end{align}
The quadratic function $h_{\alp,\eta}(\rho)$ ($-1\le\rho\le1$) satisfies
\begin{align*}
h_{\alp,\eta}(-1)&=Q_{1,2,3;\alp}(1+\eta)^2\ge0,\\
h_{\alp,\eta}(1)&=Q_{1,2,3;\alp}(1-\eta)^2\ge0,\\
h_{\alp,\eta}(0)&=Q_{1,2,3;\alp}(1+\eta^2)-1.
\end{align*}
It is easily seen that if $\alp\le\frac{2}{3}$ then $Q_{1,2,3;\alp}\le\f12$, so $h_{\alp,\eta}(0)\le0$ for all $\eta\in[-1,1]$. Thus if $\alp\le\frac{2}{3}$ then $h_{\alp,\eta}(\rho)$ must have one root in $[-1,0]$ and one root in $[0,1]$.  The two roots are given by
\begin{align}\label{roots}
\ts\^\rho_{\alp,\eta}^{\pm}=Q_{1,2,3;\alp}\eta\pm\sqrt{(1-Q_{1,2,3;\alp})(1-Q_{1,2,3;\alp}\eta^2)};
\end{align}
note that $\^\rho_{\alp,-\eta}^{\pm}=-\^\rho_{\alp,\eta}^{\mp}$.

\begin{table}[h!]
\centering
\begin{center}
\begin{tabular}{ |c|c|c|c|c|c|c|c|c| }
%\begin{tabular}{ |p{2.4cm}||p{.6cm}|p{.6cm}|p{.6cm}|p{.6cm}|p{.6cm}|p{.6cm}|p{.6cm}|p{.6cm}|p{.6cm}|p{.6cm}|p{.6cm}| }
% \hline
%q&\|a^{(qœ§ )}\|&\^y^{(q)}'&y^{(q)}'&3&4&5&6&7&8&9&10\\
 \hline
 $\alp$&$Q_{1,2,3;\alp}$&$\eta$&$(\^\rho_{\alp,\eta}^-,\,\^\rho_{\alp,\eta}^+)$\\
 \hline
 \hline
 .50&.333&  1&(-.333,\,1)\\
  &  &.5&(-.615,\,.948)\\
  &  &0&(-.816,\,.816)\\
    \hline
.20&.111&1&(-.778,\,1)\\
&&.5&(-.874,\,.985)\\
&&0&(-.943,\,.943)\\
 \hline
.10&.0526&1&(-.895,\,1)\\
&&.5&(-.941,\,.993)\\
&&0&(-.973,\,.973)\\
 \hline
.05&.0256&1&(-.949,\,1)\\
&&.5&(-.971,\,.997)\\
&&0&(-.987,\,.987)\\
 \hline
.01&.00503&1&(-.990,\,1)\\
&&.5&(-.994,\,.999)\\
&&0&(-.997,\,.997)\\
 \hline
 0+&0&1&(-1,\,1)\\
&&.5&(-1,\,1)\\
&&0&(-1,\,1)\\
 \hline
\end{tabular}
\end{center}
\caption{For $p=2$, $N=3$, and sufficiently large $\gam_1^2\equiv\max(\Lam_1,\Lam_2)$ (i.e., $\gam_1^2\ge\frac{1}{3}\Lam_{1,2,3;\alp}^*$), if $\rho\in(\^\rho_{\alp,\eta}^-,\,\^\rho_{\alp,\eta}^+)$ then $|\^\ome_\alp(\gam_1,\eta,\rho)|=1$, i.e., the power of at least one of the two univariate size-$\alp$ $t^2$-tests dominates that of the bivariate size-$\alp$  $T^2$-test. (Note that $\^\rho_{\alp,-\eta}^{\pm}=-\^\rho_{\alp,\eta}^{\mp}$.)}
\label{table:1}
\end{table}

It follows that if $\alp\le\frac{2}{3}$ then for sufficiently large $\gam_1^2$, i.e., $\gam_1^2\ge\frac{1}{3}\Lam_{1,2,3;\alp}^*$,
\begin{equation}\label{asympt8}
\ts\rho\in(\^\rho_{\alp,\eta}^-,\,\^\rho_{\alp,\eta}^+)\implies |\^\ome_\alp(\gam_1,\eta,\rho)|=1,
%\pi_\alp(\gam^2;1,2)>\pi_\alp\big(\ts\bet_\alp\big(\frac{1-2\eta\rho+\eta^2}{1-\rho^2}\big);2,1\big).
\end{equation}
that is, at least one of the two univariate $t^2$-tests is more powerful than the overall (bivariate) $T^2$-test. Specifically, when $\gam_1^2\ge\frac{1}{3}\Lam_{1,2,3;\alp}^*$, $|\^\ome_\alp(\gam_1,\eta,\rho)|=1$ in the $(\eta,\rho)$-regions of the parameter space indicated in Table 1. From this it is seen that for $p=2$, $N=3$, large $\gam_1^2$, and the common (small) values of $\alp$, the overall size-$\alp$ $T^2$-test is dominated by at least one of the two univariate size-$\alp$  $t^2$-tests over a large portion of the alternative hypothesis $K$.\hfill$\square$
\vskip4pt

\newpage

\nid{\bf Example 2.2.} Suppose that $p\ge3$ and $N=p+2$.
%Abbreviate the 3 non-empty subsets $\{1\}$, $\{2\}$, and $\{12\}$ of $\Ome_2$ by $1$, $2$, and $12$ respectively. 
The powers of the ${p\choose2}$ bivariate size-$\alp$ $T^2$-tests and the overall $p$-variate size-$\alp$ $T^2$-test will be compared for distant alternatives, which requires comparison of the powers
\begin{align}
\{\pi(\Lam_\ome;2,p)\mid\ome\in\Ome_p,\,|\ome|=2\}\ \  \mathrm{and}\ \ \pi(\Lam;p,\,2).\label{powers}
\end{align}
From \eqref{asymp3},
\begin{align*}
\Lam^{(2)}:&=\max\{\Lam_\ome\mid\ome\in\Ome_p,\,|\ome|=2\}>
\max\left(\Lam_{2,p,p+2;\alp}^*,\;Q_{2,p,p+2;\alp}\Lam\right)\\%\label{asymp3a}\\
&\ \ \ \ \implies\ \max\{\pi(\Lam_\ome;2,p)\mid\ome\in\Ome_p,\,|\ome|=2\}>\pi(\Lam;p,\,2).%\label{asymp3b}
\end{align*}
Therefore for sufficiently large values of $\Lam^{(2)}$, namely $\Lam^{(2)}\ge\Lam_{2,p,p+2;\alp}^*$, at least one of the bivariate size-$\alp$ $T^2$-tests will be more powerful than the overall (p-variate) size-$\alp$ $T^2$-test\footnote{\label{Foot20} Note that by itself this does not establish that $|\^\ome_\alp(\gam,R)|=2$.} provided that
\begin{equation}\label{Lam2greater}
\ts\Lam^{(2)}>Q_{2,p,p+2;\alp}\Lam.%\left(\frac{1-(1-\alp)^{\frac{2}{N-2}}}{\alp^{\frac{2}{N-2}}}\right)\lam_{N-2}.
\end{equation}

From \eqref{twobs} we obtain the explicit expression
\begin{equation}\label{Q12A}
\ts Q_{2,p,p+2;\alp}=\frac{b_{2,p;\alp}}{b_{p,2;\alp}}=\frac{1-(1-\alp)^{\frac{2}{p}}}{\alp^{\frac{2}{p}}}.
\end{equation}
If we set $\nu_p=\frac{2}{p}$ and $U_{p;\alp}=\frac{\nu_p}{Q_{2,p,p+2;\alp}}$ then
\begin{align}
\lim\limits_{p\to\infty}Q_{2,p,p+2;\alp}&=0,\label{limnuQ1}\\
\lim_{p\to\infty}U_{p;\alp}&\ts=\frac{1}{-\log(1-\alp)}\ \ \ (>1\ \ \mathrm{for}\ \ \alp<\frac{e-1}{e}=.6321);\label{limnuQ}\\
\lim_{p\to\infty}(p-1)U_{p;\alp}&=\infty.\label{limU}
\end{align}
Table 2 shows that $Q_{2,p,p+2;\alp}$ decreases rapidly to 0 as $p\to\infty$, which suggests that \eqref{Lam2greater} might hold over substantial regions of the alternative hypothesis $K$. We proceed to exhibit several such regions.
\vskip4pt

\begin{table}[h!]
\centering
\begin{center}
\begin{tabular}{ |c|c|c|c|c|c|c|c|c| }
%\begin{tabular}{ |p{2.4cm}||p{.6cm}|p{.6cm}|p{.6cm}|p{.6cm}|p{.6cm}|p{.6cm}|p{.6cm}|p{.6cm}|p{.6cm}|p{.6cm}|p{.6cm}| }
% \hline
%q&\|a^{(qœ§ )}\|&\^y^{(q)}'&y^{(q)}'&3&4&5&6&7&8&9&10\\
 \hline
 $p$&$\alp$&$Q_{2,p,p+2;\alp}$&$U_{p;\alp}$&$\~\psi_{p;\alp}^-$&$(\~\rho_{p;\alp}^-,\,\~\rho_{p;\alp}^+)$&$-\frac{1}{p-1}$\\
 \hline
 \hline
  4&  .20&.236&2.118&-.209&(-.304,\,.776)&-.333\\
 &  .10&.162&3.081&-.245&(-.314,\,.844)&\\
  &  .05&.113&4.415&-.279&(-.320,\,.890)&\\
  &  .01&.050&9.975&-.310&(-.328,\,.951)&\\
    \hline
 10&.20&.0602&3.321&-.080&(-.110,\,.907)&-.111\\
&.10&.0330&6.052&-.094&(-.110,\,.948)&\\
&.05&.0186&10.764&-.102&(-.111,\,.973)&\\
&.01&$.0^2504$&39.652&-.109&(-.111,\,.980)&\\
 \hline
 20&.20&.0259&3.858&-.039&(-.0525,\,.954)&-.0526\\
&.10&.0132&7.579&-.046&(-.0526,\,.976)&\\
&.05&$.0^2690$&14.486&-.049&(-.0526,\,.988)&\\
&.01&$.0^2159$&62.810&-.052&(-.0526,\,.997)&\\
 \hline
40&.20&.0120&4.158&-.020&(-.0256,\,.977)&-.0256\\
&.10&$.0^2590$&8.481&-.023&(-.0256,\,.989)&\\
&.05&$.0^2298$&16.806&-.024&(-.0256,\,.994)&\\
&.01&$.0^3632$&79.055&-.025&(-.0256,\,.999)&\\

 \hline
 $\infty$&.20&0&4.481&0&[0,\,1)&0\\
&.10&0&9.491&0&[0,\,1)&\\
&.05&0&19.496&0&[0,\,1)&\\
&.01&0&99.499&0&[0,\,1)&\\
 \hline
\end{tabular}
\end{center}
\caption{Take $p\ge3$, $N=p+2$, and $\del^2$ sufficiently large, that is, $\del^2\ge (p+2)^{-1}\Lam_{2,p,p+2;\alp}^*$. In Case 1 (Case 2) if $\~\psi_{p;\alp}^-<\!\rho\!<1$ ($\~\rho_{p;\alp}^-<\rho<\~\rho_{p;\alp}^+$) then the power of at least one of the ${p\choose2}$ bivariate size-$\alp$ $T^2$-tests dominates that of the overall size-$\alp$ $T^2$-test. Note that the feasible range of $\rho$ is $\big(-\frac{1}{p-1},1\big)$.}
\label{table:2}
\end{table}

\nid{\it Case 1: $\gam_1=\cdots=\gam_p=:\del$ and $R$ has the intraclass form}
\begin{equation}\label{intraclass}
\ts R_\rho:=(1-\rho)I_p+\rho {\bf e}_p{\bf e}_p'.
\end{equation}
Here $-\frac{1}{p-1}<\rho<1$ and ${\bf e}_p=(1,\dots,1)':p\times1$.  Then
\begin{align}
R_\rho^{-1}&\ts=\frac{1}{1-\rho}\left[I_p-\frac{\rho {\bf e}_p{\bf e}_p'}{1+\rho(p-1)}\right],\label{Rinverse}\\
{\bf e}_p'R_\rho^{-1}{\bf e}_p&\ts=\frac{p}{1+\rho(p-1)},\label{eRinve}
\end{align}
so by \eqref{lamgamR} and \eqref{lamome},
\begin{align}
\ts\Lam^{(2)}=\frac{2(p+2)\del^2}{1+\rho},\qquad\Lam=\frac{p(p+2)\del^2}{1+(p-1)\rho}.\label{Lam2lam}
\end{align}
Thus $\Lam^{(2)}\ge\Lam_{2,p,p+2;\alp}^*$ holds for all feasible $\rho$ if $\del^2\ge (p+2)^{-1}\Lam_{2,p,p+2;\alp}^*$.
Also, if we set $\nu_p=\frac{2}{p}\ (\le\frac{2}{3})$ then \eqref{Lam2greater} is equivalent to each of the inequalities
\begin{align}
\ts U_{p;\alp}&>\ts\frac{1+\rho}{1+(p-1)\rho};\nonumber\\%\label{Lam2greater1}\\
\ts\left[(p-1)U_{p;\alp}-1\right]\rho&>1-U_{p;\alp}.\label{Lam2greater2}
\end{align}
Because $(p-1)U_{p;\alp}>1$ for common (small) values of $\alp$ (see \eqref{limU} and Table 2), in such cases \eqref{Lam2greater2} in turn is equivalent to
\begin{align}
\rho&>\ts\frac{1-U_{p;\alp}}{(p-1)U_{p;\alp}-1}=:\~\psi_{p;\alp}^-.\label{Lam2greater3}
\end{align}
Table 2 shows that in Case 1, $\~\psi_{p;\alp}^-$ is close to the lower limit of the feasible range $(-\frac{1}{p-1},1)$ for $\rho$. Thus by \eqref{Lam2greater3}, if $\del^2\ge (p+2)^{-1}\Lam_{2,p,p+2;\alp}^*$ then at least one of the bivariate size-$\alp$ $T^2$-tests will be more powerful than the overall (p-variate) size-$\alp$ $T^2$-test for most of the region in the alternative hypothesis specified in Case 1.
\vskip4pt

\newpage

\nid{\it Case 2: $\gam_i=\gam_j=:\del$ for some $\{i,j\}\subset\{1,\dots,p\}$, $\gam_k=0$ for $k\ne i,j$, and $R$ has the intraclass form $R_\rho$ in \eqref{intraclass}.} By \eqref{lamgamR} and \eqref{lamome},
\begin{align}
\ts\Lam^{(2)}=\frac{2(p+2)\del^2}{1+\rho},\qquad\Lam=\frac{2(p+2)\del^2[1+(p-3)\rho]}{(1-\rho)[1+(p-1)\rho]},\label{Lam2lamA}
\end{align}
so again $\Lam^{(2)}\ge\Lam_{2,p,p+2;\alp}^*$ holds for all feasible $\rho$ if $\del^2\ge (p+2)^{-1}\Lam_{2,p,p+2;\alp}^*$. Also, abbreviating $Q_{2,p,p+2;\alp}$ by $Q$, \eqref{Lam2greater} is equivalent to each of the inequalities
\begin{align*}
\ts\frac{1}{1+\rho}&>%\left(\frac{1-(1-\alp)^{\nu_N}}{\alp^{\nu_N}}\right)
\ts\frac{Q[1+(p-3)\rho]}{(1-\rho)[1+(p-1)\rho]};\\%\nonumber\\%\label{Lam2greater4}\\
%0&>(1+\rho)[1+(N-3)\rho]-\ts \frac{\nu_N}{Q_{2,N-2,N;\alp}},\\
 0&>\ts[(p-1)+(p-3)Q]\rho^2-(p-2)(1-Q)\rho-(1-Q)=:h_{p;\alp}(\rho).
 %\label{Lam2greater5}
\end{align*}
Since  $h_{p;\alp}(0)=Q-1<0$ for common (small) values of $\alp$ (see \eqref{limnuQ1}-\eqref{limnuQ} and Table 2), $h_{p;\alp}(\rho)$ has two real roots $\~\rho_{p;\alp}^-<0<\~\rho_{p;\alp}^+$ (found numerically). Therefore $0>h_{p;\alp}(\rho)$ for $\~\rho_{p;\alp}^-<\rho<\~\rho_{p;\alp}^+$. 

Table 2 shows that in Case 2, the interval $(\~\rho_{p;\alp}^-,\~\rho_{p;\alp}^+)$ covers almost all of the feasible range $(-\frac{1}{p-1},1)$ for $\rho$. Thus 
%by \eqref{Lam2greater5}, 
if $\del^2\ge (p+2)^{-1}\Lam_{2,p,p+2;\alp}^*$ then at least one of the bivariate size-$\alp$ $T^2$-tests will be more powerful than the overall p-variate size-$\alp$ $T^2$-test for most of the region in the alternative hypothesis specified by Case 2.

\vskip4pt

\nid{\it Case 3: $\gam_i=\del$ and $\gam_j=-\del$ for some $\{i,j\}\subset\{1,\dots,p\}$, $\gam_k=0$ for $k\ne i,j$, and $R$ has the intraclass form $R_\rho$.} By \eqref{lamgamR} and \eqref{lamome},
\begin{align}
\ts\Lam^{(2)}=\frac{2(p+2)\del^2}{1-\rho},\qquad\Lam=\frac{2(p+2)\del^2}{1-\rho}.\label{Lam2lamB}
\end{align}
Thus $\Lam^{(2)}\ge\Lam_{2,p,p+2;\alp}^*$ again holds for all feasible $\rho$ if $\del^2\ge (p+2)^{-1}\Lam_{2,p,p+2;\alp}^*$, while \eqref{Lam2greater} is equivalent to $1>Q_{2,p,p+2;\alp}$, which holds for most $p,\alp$ (see  \eqref{limnuQ1}-\eqref{limnuQ} and Table 2). Therefore if $\del^2\ge (p+2)^{-1}\Lam_{2,p,p+2;\alp}^*$ then at least one of the bivariate size-$\alp$ $T^2$-tests will be more powerful than the overall p-variate size-$\alp$ $T^2$-test for the {\it entire} region in the alternative hypothesis covered by Case 3.
\vskip4pt

\nid{\it Case 4: $p=:2l$ is even, $\gam_i=\del$ for $l$ indices in $\{1,\dots,p\}$, $\gam_i=-\del$ for the remaining $l$ indices, and $R$ has the intraclass form $R_\rho$.} By \eqref{lamgamR} and \eqref{lamome},
\begin{align}
\ts\Lam^{(2)}=\frac{2(p+2)\del^2}{1-|\rho|},\qquad\Lam=\frac{p(p+2)\del^2}{(1-\rho)}.\label{Lam2lamC}
\end{align}
Thus $\Lam^{(2)}\ge\Lam_{2,p,p+2;\alp}^*$ again holds for all feasible $\rho$ if $\del^2\ge (p+2)^{-1}\Lam_{2-2,p+2;\alp}^*$, while \eqref{Lam2greater} is equivalent to the inequality
\begin{align}
U_{p;\alp}&>\ts\frac{1-|\rho|}{1-\rho}.\label{Lam2greater6}
\end{align}
Because $\frac{1-|\rho|}{1-\rho}\le1$, while $U_{p;\alp}>1$ for holds for most $p,\alp$ (see \eqref{limnuQ} and Table 2), we see that if $\del^2\ge (p+2)^{-1}\Lam_{2,p,p+2;\alp}^*$ then  at least one of the bivariate size-$\alp$ $T^2$-tests is more powerful than the overall p-variate size-$\alp$ $T^2$-test over the {\it entire} region in the alternative hypothesis in Case 4.\hfill$\square$
\vskip4pt

%COMMENT
\begin{comment}
\nid{\bf OMIT-too Similar to Case 2}{\it Case 3: $\gam_1=\del$, $\gam_2^2=\cdots=\gam_{N-2}^2=0$.} By \eqref{lamgamR} and \eqref{lamome},
\begin{align*}
\ts\Lam^{(2)}=N\del^2,\qquad\Lam=\frac{N\del^2[1+(N-4)\rho]}{(1-\rho)[1+(N-3)\rho]}.
\end{align*}
\begin{align}
1&>%\left(\frac{1-(1-\alp)^{\nu_N}}{\alp^{\nu_N}}\right)
\ts\frac{Q[1+(N-4)\rho]}{(1-\rho)[1+(N-3)\rho]};\nonumber\\%\label{Lam2greater4}\\
%0&>(1+\rho)[1+(N-3)\rho]-\ts \frac{\nu_N}{Q_{2,N-2,N;\alp}},\\
 0>\ts[(N-3)&]\rho^2-(N-4)(1-Q)\rho-(1-Q)=:h_{N;\alp}(\rho).\label{Lam2greater5}
\end{align}
\end{comment}
%END COMMENT

\nid{\bf 3. Some local power comparisons.} From \eqref{pialp}-\eqref{cmnk} and \eqref{cmnk2}, as $\Lam_\ome\downarrow0$ the power function $\pi_\alp(\Lam_\ome)\equiv\pi_\alp(\Lam_\ome;|\ome|,N-|\ome|)$ of the $T_\ome^2$-test satisfies 
\begin{align}
\pi_\alp(\Lam_\ome)&\ts=e^{-\frac{\Lam_\ome}{2}}[\alp+\frac{\Lam_\ome}{2}c_{|\ome|,N-|\ome|;1;\alp}+O(\Lam_\ome^2)]\label{local1}\\
 &=\ts\alp+\frac{\Lam_\ome}{2}(c_{|\ome|,N-|\ome|;1;\alp}-\alp)+O(\Lam_\ome^2).\label{local12}
\end{align}
Thus for two subsets $\ome,\,\ome'$ s.t. $\ome\subset\ome'$,
%$|\ome|\subset|\ome'|$,  
$\exists\;\Lam_{|\ome|,|\ome'|,N;\alp}^{**}>0$ s.t.
%$\exists\;\Lam_\alp\equiv\Lam_{|\ome|,|\ome'|;\alp}>0$ s.t.
\begin{equation}\label{asymp35}
\ts\Lam_{\ome'}<\min\left(\Lam_{|\ome|,|\ome'|,N;\alp}^{**},\;Z_{|\ome|,|\ome'|,N;\alp}\Lam_{\ome}\right)
\ \implies\ \pi_\alp(\Lam_\ome)>\pi_\alp(\Lam_{\ome'}),
\end{equation}
where, from (2.2) and (2.3) in [DGP],
\begin{equation}\label{Zdef}
\ts Z_{|\ome|,|\ome'|,N;\alp}:=\frac{c_{|\ome|,N-|\ome|;1;\alp}-\alp}{c_{|\ome'|,N-|\ome'|;1;\alp}-\alp}>1.
\end{equation}
Therefore power comparisons of $T_\ome^2$ and $T_{\ome'}^2$ for local alternatives\footnote{\label{Foot23} We do not claim to know the values of $\Lam_{|\ome|,|\ome'|,N;\alp}^{**}$, even approximately.} require determination of the lower tail probabilities $c_{m,n;k;\alp}$ (see \eqref{cmnk2}), which in turn require the lower quantiles $b_{n,m;\alp}$. 

In parallel with Section 2, this is done explicitly in Examples 3.1 and 3.2. As in Examples 2.1 and 2.2 these examples begin to suggest that variable selection might be limited to very small subsets $\ome\in\Ome_p$, e.g., singletons in the bivariate Example 3.1, or pairs (plus singletons) in Example 3.2.
\vskip4pt

\nid{\bf Example 3.1.} As in Example 2.1 consider the bivariate case $p=2$. Repeat the first two paragraphs from Example 2.1 verbatim, except replace ``distant alternatives" by ``local alternatives". Because $\max(\Lam_1,\Lam_2)\le\Lam$, it follows from \eqref{lams} and \eqref{asymp35} that
%for $j=1,2$,
\begin{align}
&\ts\qquad\big(\frac{1-2\eta\rho+\eta^2}{1-\rho^2}\big)\gam_1^2<\min\left(\frac{1}{N}\Lam_{1,2,N;\alp}^{**},\;Z_{1,2,N;\alp}\gam_1^2\right)\label{asymp44}\\
\implies\ &\max(\pi_\alp(\Lam_1;1,2),\pi_\alp(\Lam_2;1,2))>\pi_\alp(\Lam;2,1),\label{asymp4a4}\\
\implies\ & |\^\ome_\alp(\gam_1,\eta,\rho)|=1.\label{hatome2}
\end{align}

In the simplest case $N=3$, it follows from \eqref{cmnk2}, \eqref{betapdf}, and \eqref{twobs} that
\begin{align}
Z_{1,2,3;\alp}&\ts=\frac{c_{1,2;1;\alp}-\alp}{c_{2,1;1;\alp}-\alp}=\frac{1-(1-\alp)^3-\alp}{\frac{3}{2}(\alp-\frac{1}{3}\alp^3)-\alp}=\frac{2(2-\alp)}{1+\alp}.\label{Zdef123}
%&\ts=\frac{1-(1-\alp)^3-\alp}{\frac{3}{2}(\alp-\frac{1}{3}\alp^3)-\alp}\\
%&\ts=\frac{2(2-\alp)}{1+\alp}\\
\end{align}
First note that in \eqref{asymp44},
\begin{align}
&\ts\frac{1-2\eta\rho+\eta^2}{1-\rho^2}<Z_{1,2,3;\alp}\label{lamalp33}\\
\iff h_{\alp,\eta}(\rho):=&Z_{1,2,3;\alp}\rho^2-2\eta\rho+\eta^2+1-Z_{1,2,3;\alp}<0.\label{lamalp3}
\end{align}
The quadratic function $h_{\alp,\eta}(\rho)$ ($-1\le\rho\le1$) satisfies
\begin{align*}
h_{\alp,\eta}(-1)&=(1+\eta)^2\ge0,\\
h_{\alp,\eta}(1)&=(1-\eta)^2\ge0,\\
h_{\alp,\eta}(0)&=\eta^2+1-Z_{1,2,3;\alp},
\end{align*}
It is easily seen that if $\alp\le\frac{1}{2}$ then $Z_{1,2,3;\alp}\ge2$, so $h_{\alp,\eta}(0)\le0$ for all $\eta\in[-1,1]$. Therefore, if $\alp\le\frac{1}{2}$ then $h_{\alp,\eta}(\rho)$ must have one root in $[-1,0]$ and one root in $[0,1]$.  The two roots are given by
\begin{align}\label{rootsA}
\ts\c\rho_{\alp,\eta}^{\pm}=\frac{\eta}{Z_{1,2,3;\alp}}\pm\frac{1}{Z_{1,2,3;\alp}}\sqrt{(Z_{1,2,3;\alp}-1)(Z_{1,2,3;\alp}-\eta^2)};
\end{align}
again $\c\rho_{\alp,-\eta}^{\pm}=-\c\rho_{\alp,\eta}^{\mp}$.
Thus, if $\alp\le\frac{1}{2}$ and $\rho\in(\c\rho_{\alp,\eta}^-,\,\c\rho_{\alp,\eta}^+)$ then \eqref{lamalp33} must hold.  

To conclude that $|\^\ome_\alp(\gam_1,\eta,\rho)|=1$, $\gam_1^2$ must be sufficiently small, i.e., 
\begin{equation}\label{gamlamstar3}
\ts0<\gam_1^2<\big(\frac{1-\rho^2}{1-2\eta\rho+\eta^2}\big)\frac{\Lam_{1,2,3;\alp}^{**}}{3}.
\end{equation}
Because
\begin{equation*}
\ts\min\limits_{|\eta|\le1}\big(\frac{1-\rho^2}{1-2\eta\rho+\eta^2}\big)=\t12(1-|\rho|),
\end{equation*}
for fixed $\eta$, \eqref{gamlamstar3} will be satisfied provided that 
\begin{align}
\rho&\in(\c\rho_{\alp,\eta}^-,\,\c\rho_{\alp,\eta}^+),\label{rhointerval}\\
\c m_{\alp,\eta}:&=\max(|\c\rho_{\alp,\eta}^-|,\,|\c\rho_{\alp,\eta}^+|)<1,\label{maxrhohat}\\
\gam_1^2&\ts<\frac{1}{6}(1-m_{\alp,\eta})\Lam_{1,2,3;\alp}^{**}.\label{gamlamstar4}
\end{align}
It is straightforward to show that \eqref{maxrhohat} holds for $|\eta|<1$ but not for $|\eta|=1$.

Thus, if $\alp\le\frac{1}{2}$, $|\eta|<1$, and \eqref{rhointerval}-\eqref{gamlamstar4} are satisfied then $|\^\ome_\alp(\gam_1,\eta,\rho)|=1$, in which case at least one of the two univariate $t^2$-tests are more powerful than the overall (bivariate) $T^2$-test. This occurs in the $(\eta,\rho)$-regions of the parameter space indicated in Table 3, provided that $\gam_1^2\ts<\frac{1}{6}(1-\c m_{\alp,\eta})\Lam_{1,2,3;\alp}^{**}$. Thus, for $p=2$, $N=3$, sufficiently small $\gam_1^2$, and the common (small) values of $\alp$, the overall (bivariate) size-$\alp$ $T^2$-test is dominated by at least one of the two univariate size-$\alp$  $t^2$-tests over a large portion of the alternative hypothesis.\hfill$\square$
\vskip4pt

\begin{table}[h!]
\centering
\begin{center}
\begin{tabular}{ |c|c|c|c|c|c|c|c|c| }
%\begin{tabular}{ |p{2.4cm}||p{.6cm}|p{.6cm}|p{.6cm}|p{.6cm}|p{.6cm}|p{.6cm}|p{.6cm}|p{.6cm}|p{.6cm}|p{.6cm}|p{.6cm}| }
% \hline
%q&\|a^{(qœ§ )}\|&\^y^{(q)}'&y^{(q)}'&3&4&5&6&7&8&9&10\\
 \hline
 $\alp$&$Z_{1,2,3;\alp}$&$\eta$&$(\c\rho_{\alp,\eta}^-,\,\c\rho_{\alp,\eta}^+)$&$1-\c m_{\alp,\eta}$\\
 \hline
 \hline
 .50&2&  .9&(-.095,\,.9954)&.0046\\
  &&  .5&(-.411,\,.911)&.089\\
  &&  0&(-.707,\,.707)&.293\\
    \hline
.20&3&.9&(-.398,\,.9976)&.0024\\
&&.5&(-.615,\,.948)&.052\\
&&0&(-.816,\,.816)&.184\\
 \hline
.10&3.454&.9&(-.477,\,.9980)&.0020\\
&&.5&(-.667,\,.957)&.043\\
&&0&(-.843,\,.843)&.157\\
 \hline
.05&3.714&.9&(-.514,\,.9982)&.0018\\
&&.5&(-.691,\,.962)&.038\\
&&0&(-.855,\,.855)&.145\\
 \hline
.01&3.941&.9&(-.542,\,.9983)&.0017\\
&&.5&(-.709,\,.963)&.037\\
&&0&(-.864,\,.864)&.136\\
 \hline
 0+&4&.9&(-.548,\,.9984)&.0016\\
&&.5&(-.714,\,.964)&.036\\
&&0&(-.866,\,.866)&.134\\
 \hline
\end{tabular}
\end{center}
\caption{For $p=2$, $N=3$, and sufficiently small $\gam_1^2\equiv\max(\Lam_1,\Lam_2)$ (i.e., $\gam_1^2\ts<\frac{1}{6}(1-\c m_{\alp,\eta})\Lam_{1,2,3;\alp}^{**}$), if $\rho\in(\c\rho_{\alp,\eta}^-,\,\c\rho_{\alp,\eta}^+)$ then $|\^\ome_\alp(\gam_1,\eta,\rho)|=1$, i.e., the power of at least one of the two univariate size-$\alp$ $t^2$-tests dominates that of the bivariate size-$\alp$  $T^2$-test. (Note that $\c\rho_{\alp,-\eta}^{\pm}=-\c\rho_{\alp,\eta}^{\mp}$.)}
\label{table:3}
\end{table}

\newpage

\nid{\bf Example 3.2.} Suppose that $p\ge3$ and $N=p+2\ge3$. We shall compare the powers of the ${p\choose2}$ bivariate size-$\alp$ $T^2$-tests and the overall $p$-variate size-$\alp$ $T^2$-test for local alternatives, which requires comparison of the powers
\begin{align}
\{\pi(\Lam_\ome;2,p)\mid\ome\in\Ome_p,\,|\ome|=2\}\ \  \mathrm{and}\ \ \pi(\Lam;p,\,2).\label{powersA}
\end{align}
From \eqref{asymp35},
\begin{align}
\Lam&<\min(\Lam_{2,p,p+2;\alp}^{**},\,Z_{2,p,p+2;\alp}\Lam^{(2)})\label{lamimplies2}\\
&\ \ \ \ \implies\ \max\{\pi(\Lam_\ome;2,p)\mid\ome\in\Ome_p,\,|\ome|=2\}>\pi(\Lam;p,\,2).\label{lamimplies43}
\end{align}
Therefore for sufficiently small values of $\Lam$, namely $\Lam\le\Lam_{2,p,p+2;\alp}^{**}$, at least one of the bivariate size-$\alp$ $T^2$-tests will be more powerful than the overall $p$-variate size-$\alp$ $T^2$-test\footnote{\label{Foot30} As in Example 2.2, this does not establish that $|\^\ome_\alp(\gam,R)|=2$.} whenever
\begin{equation}\label{ZLam2greater}
\Lam<Z_{2,p,p+2;\alp}\Lam^{(2)}.%\left(\frac{1-(1-\alp)^{\frac{2}{N-2}}}{\alp^{\frac{2}{N-2}}}\right)\lam_{N-2}.
\end{equation}

From \eqref{Zdef}, \eqref{cmnk2}-\eqref{betapdf}, \eqref{twobs}, and some algebra, the explicit expression
\begin{align}
Z_{2,p,p+2;\alp}&\ts\equiv\frac{c_{2,p;1;\alp}-\alp}{c_{p,2;1;\alp}-\alp}
=\ts\frac{p\alp[1-\alp^{\frac{2}{p}}]} {2(1-\alp)[1-(1-\alp)^{\frac{2}{p}}]}\label{Zdef3}
\end{align}
is obtained. Setting $\nu_p=\frac{2}{p}\ (\le\frac{2}{3})$ and $V_{p;\alp}:=\nu_p Z_{2,p,p+2;\alp}$, we have
\begin{align}
\lim\limits_{p\to\infty}Z_{2,p,p+2;\alp}&=\infty,\label{limnuZ1}\\
\lim_{p\to\infty}V_{p;\alp}&\ts=\frac{\alp\log\alp}{(1-\alp)\log(1-\alp)}\ \ \ (>1\ \ \mathrm{for}\ \ \alp<\f12),\label{limnuZ}\\
\lim_{p\to\infty}(p-1)V_{p;\alp}&=\infty.\label{limNV}
\end{align}
%\lim_{N\to\infty}(N-3)U_{N;\alp}&=\infty.\label{limU}
Table 4 shows that $Z_{2,p,p+2;\alp}$ increases rapidly to $\infty$ as $p\to\infty$, which suggests that \eqref{ZLam2greater} might hold over substantial regions of the alternative hypothesis. Several such regions are now exhibited.
\vskip4pt

\nid{\it Case 1: $\gam_1=\cdots=\gam_p=:\del$ and $R$ has the intraclass form \eqref{intraclass}.}
%\begin{equation}\label{intraclass}
%\ts R_\rho:=(1-\rho)I_{N-2}+\rho {\bf e}_{N-2}{\bf e}_{N-2}'.
%\end{equation}
Here $-\frac{1}{p-1}<\rho<1$ and as in \eqref{Lam2lam},
\begin{align}
\ts\Lam^{(2)}=\frac{2(p+2)\del^2}{1+\rho},\qquad\Lam=\frac{p(p+2)\del^2}{1+(p-1)\rho}.\label{Lamlamrepeat}
\end{align}
Here \eqref{ZLam2greater} is equivalent to each of the inequalities
\begin{align}
\ts V_{p;\alp}&>\ts\frac{1+\rho}{1+(p-1)\rho};\nonumber\\%\label{Lam2greater1}\\
\ts V_{p;\alp}-1&>-\left[(p-1)V_{p;\alp}-1\right]\rho.\label{ZLam2greater2}
\end{align}
Because $(p-1)V_{p;\alp}>1$ for common (small) values of $\alp$ (see \eqref{limNV} and Table 4), in such cases \eqref{ZLam2greater2} in turn is equivalent to
\begin{align}
\u\psi_{p;\alp}^-:=\ts-\frac{V_{p;\alp}-1}{(p-1)V_{p;\alp}-1}<\rho<-1.\label{ZLam2greater3}
\end{align}

To conclude that \eqref{lamimplies2}-\eqref{lamimplies43} holds, $\del^2$ must be sufficiently small, i.e.,
\begin{equation}\label{del2lele}
\ts0<\del^2< \Big[\frac{1+(p-1)\rho}{p(p+2)}\Big]\Lam_{2,p,p+2;\alp}^{**}.
\end{equation}
However, $\rho>\u\psi_{p;\alp}^-$ implies that
\begin{equation}\label{ineq10}
\ts\frac{1+(p-1)\rho}{p(p+2)}>\frac{1+(p-1)\u\psi_{p;\alp}^-}{p(p+2)}
=\frac{p-2}{p(p+2)[(p-1)V_{p;\alp}-1]}:=\u m_{p,\alp}.
\end{equation}
Therefore \eqref{del2lele} will be satisfied provided that 
\begin{align}
\rho&>\u\psi_{p;\alp}^-\quad\mathrm{and}\quad\del^2\ts<\u m_{p,\alp}\Lam_{2,p,p+2;\alp}^{**}.\label{rhodel}
%\del^2&\ts<m_{N,\alp}\lam_{2,N-2,N;\alp}^{**}.\label{gamlamstar4}
\end{align}

If $p$ is large and $\alp$ is small, Table 4 shows that in Case 1, $\u\psi_{p;\alp}^-$ is close to the lower limit of the feasible range $(-\frac{1}{p-1},1)$ for $\rho$. Here, by \eqref{rhodel}, if $\del^2\ts<\u m_{p,\alp}\Lam_{2,p,p+2;\alp}^{**}$ then at least one of the bivariate size-$\alp$ $T^2$-tests will be more powerful than the overall ($p$-variate) size-$\alp$ $T^2$-test for most of the region in the alternative hypothesis covered by Case 1.

\vskip4pt

\nid{\it Case 2: $\gam_i=\gam_j=:\del$ for some $\{i,j\}\subset\{1,\dots,p\}$, $\gam_k=0$ for $k\ne i,j$, and $R$ has the intraclass form $R_\rho$ in \eqref{intraclass}.} As in \eqref{Lam2lamA},
\begin{align*}
\ts\Lam^{(2)}=\frac{2(p+2)\del^2}{1+\rho},\qquad\Lam=\frac{2(p+2)\del^2[1+(p-3)\rho]}{(1-\rho)[1+(p-1)\rho]},
\end{align*}
Abbreviating $Z_{2,p,p+2;\alp}$ by $Z$, \eqref{ZLam2greater} is equivalent to each of the inequalities
\begin{align*}
\ts\frac{Z}{1+\rho}&>\ts\frac{[1+(p-3)\rho]}{(1-\rho)[1+(p-1)\rho]};\\
 0&>\ts[(p-1)Z+(p-3)]\rho^2-(p-2)(Z-1)\rho-(Z-1)=:h_{p;\alp}(\rho).
 %\label{Lam2greater5}
\end{align*}
Since $h_{p;\alp}(0)=1-Z<0$ (cf. \eqref{Zdef}), $h_{p;\alp}(\rho)$ has real roots $\u\rho_{p;\alp}^-<0<\u\rho_{p;\alp}^+$ (found numerically). Therefore $0>h_{p;\alp}(\rho)$ for $\u\rho_{p;\alp}^-<\rho<\u\rho_{p;\alp}^+$. 

To conclude that \eqref{lamimplies2}-\eqref{lamimplies43} holds, $\del^2$ must be sufficiently small, i.e.,
\begin{equation}\label{del2lele2}
\ts0<\del^2< \Big\{\frac{(1-\rho)[1+(p-1)\rho]}{2(p+2)[1+(p-3)\rho]}\Big\}\Lam_{2,p,p+2;\alp}^{**}.
\end{equation}
Because $\frac{(1-\rho)[1+(p-1)\rho]}{1+(p-3)\rho}$ is decreasing in $\rho$, $\rho<\u\rho_{p;\alp}^+$ implies that
\begin{equation}\label{ineq10A}
\ts\frac{(1-\rho)[1+(p-1)\rho]}{2(p+2)[1+(p-3)\rho]}>\frac{(1-\u\rho_{p;\alp}^+)[1+(p-1)\u\rho_{p;\alp}^+]}{2(p+2)[1+(p-3)\u\rho_{p;\alp}^+]}:=\u m_{p,\alp}'.
\end{equation}
Therefore \eqref{del2lele2} will be satisfied provided that 
\begin{align}
\rho\in(\u\rho_{p;\alp}^-,\,\u\rho_{p;\alp}^+)\quad\mathrm{and}\quad\del^2<\u m_{p,\alp}'\Lam_{2,p,p+2;\alp}^{**}.\label{rhodel2}
%\del^2&\ts<m_{N,\alp}\lam_{2,N-2,N;\alp}^{**}.\label{gamlamstar4}
\end{align}

Table 4 shows that in Case 2, the interval $(\u\rho_{p;\alp}^-,\u\rho_{p;\alp}^+)$ covers almost all of the feasible range $(-\frac{1}{p-1},1)$ for $\rho$. Thus if $\del^2<\u m_{p,\alp}'\Lam_{2,p,p+2;\alp}^{**}$ then 
at least one of the bivariate size-$\alp$ $T^2$-tests will be more powerful than the overall p-variate size-$\alp$ $T^2$-test for most of the region in the alternative hypothesis determined by Case 2.
\vskip4pt

\nid{\it Case 3: $\gam_i=\del$ and $\gam_j=-\del$ for some $\{i,j\}\subset\{1,\dots,p\}$, $\gam_k=0$ for $k\ne i,j$, and $R$ has the intraclass form $R_\rho$.} By \eqref{lamgamR} and \eqref{lamome},
\begin{align*}
\ts\Lam^{(2)}=\frac{2(p+2)\del^2}{1-\rho},\qquad\Lam=\frac{2(p+2)\del^2}{1-\rho}.
\end{align*}
Here \eqref{ZLam2greater} is equivalent to $Z_{2,p,p+2;\alp}>1$, which holds for all $p,\alp$ (see  \eqref{Zdef}). 

To conclude that \eqref{lamimplies2}-\eqref{lamimplies43} holds, $\del^2$ must be sufficiently small, i.e.,
\begin{equation}\label{del2lele3}
\ts\del^2< \big(\frac{1-\rho}{2(p+2)}\big)\Lam_{2,p,p+2;\alp}^{**}
\end{equation}
for all $\rho\in(-\frac{1}{p-1},1)$. This requires that $\rho$ be bounded below 1, that is, $\rho<1-\eps$ for some $\eps>0$, whence \eqref{del2lele3} will be satisfied if 
\begin{equation}\label{del2eps2N}
\ts\del^2<\big(\frac{\eps}{2(p+2)}\big)\Lam_{2,p,p+2;\alp}^{**}.
\end{equation}
Therefore if \eqref{del2eps2N} holds then at least one of the bivariate size-$\alp$ $T^2$-tests will be more powerful than the overall $p$-variate size-$\alp$ $T^2$-test if $\rho<1-\eps$, which covers almost all of the region in the alternative hypothesis determined by Case 3.
\vskip4pt

\nid{\it Case 4: $p=:2l$ is even, $\gam_i=\del$ for $l$ indices in $\{1,\dots,p\}$, $\gam_i=-\del$ for the remaining $l$ indices, and $R$ has the intraclass form $R_\rho$.} By \eqref{lamgamR} and \eqref{lamome},
\begin{align*}
\ts\Lam^{(2)}=\frac{2(p+2)\del^2}{1-|\rho|},\qquad\Lam=\frac{p(p+2)\del^2}{(1-\rho)}.
\end{align*}
Here \eqref{ZLam2greater} is equivalent to the inequality
\begin{align}
V_{p;\alp}&>\ts\frac{1-|\rho|}{1-\rho}.\label{ZLam2greater6}
\end{align}
Because $\frac{1-|\rho|}{1-\rho}\le1$, while $V_{p;\alp}>1$ for holds for most $p,\alp$ (see \eqref{limnuZ} and Table 4), \eqref{ZLam2greater6} is satisfied for most $p,\alp$. 

To conclude that \eqref{lamimplies2}-\eqref{lamimplies43} holds, $\del^2$ must be sufficiently small, i.e.,
\begin{equation}\label{del2lele5}
\ts\del^2< \big[\frac{1-\rho}{p(p+2)}\big]\Lam_{2,p,p+2;\alp}^{**}
\end{equation}
for all $\rho\in(-\frac{1}{p-1},1)$. This again requires that $\rho$ be bounded below 1, that is, $\rho<1-\eps$ for some $\eps>0$, whence \eqref{del2lele5} will be satisfied if 
\begin{equation}\label{del2lele6}
\ts\del^2< \big[\frac{\eps}{p(p+2)}\big]\Lam_{2,p,p+2;\alp}^{**}
\end{equation}
Therefore if \eqref{del2lele6} holds then at least one of the bivariate size-$\alp$ $T^2$-tests will be more powerful than the overall $p$-variate size-$\alp$ $T^2$-test if $\rho<1-\eps$, which covers almost all of the region in the alternative hypothesis determined by Case 4.\hfill$\square$
\vskip4pt

\begin{table}[h!]
\centering
\begin{center}
\begin{tabular}{ |c|c|c|c|c|c|c|c|c| }
%\begin{tabular}{ |p{2.4cm}||p{.6cm}|p{.6cm}|p{.6cm}|p{.6cm}|p{.6cm}|p{.6cm}|p{.6cm}|p{.6cm}|p{.6cm}|p{.6cm}|p{.6cm}| }
% \hline
%q&\|a^{(qœ§ )}\|&\^y^{(q)}'&y^{(q)}'&3&4&5&6&7&8&9&10\\
 \hline
 $p$&$\alp$&$Z_{2,p,p+2;\alp}$&$V_{p;\alp}$&$\u\psi_{p;\alp}^-$&$\u m_{p,\alp}$&$(\u\rho_{p;\alp}^-,\u\rho_{p;\alp}^+)$&$\u m_{p,\alp}'$&$-\frac{1}{p-1}$\\
 \hline
 \hline
  4&  .20&2.618&1.309&-.106&.0285&(-.282,\,.648)&.0524&-.333\\
 &  .10&2.961&1.480&-.140&.0242&(-.289,\,.686)&.0475&\\
  &  .05&3.228&1.614&-.160&.0217&(-.295,\,.710)&.0442&\\
  &  .01&3.627&1.814&-.183&.0118&(-.298,\,.741)&.0400&\\
    \hline
 10&.20&7.882&1.576&-.044&$.0^2506$&(-.108,\,.815)&$.0^2958$&-.111\\
&.10&9.832&1.966&-.058&$.0^2399$&(-.109,\,.849)&$.0^2783$&\\
&.05&11.621&2.324&-.066&$.0^2335$&(-.109,\,.871)&$.0^2669$&\\
&.01&15.138&3.028&-.077&$.0^2254$&(-.110,\,.899)&$.0^2630$&\\
 \hline
 20&.20&16.842&1.684&-.022&$.0^2132$&(-.0523,\,.898)&$.0^2257$&-.0526\\
&.10&21.804&2.180&-.029&$.0^2101$&(-.0524,\,.921)&$.0^2199$&\\
&.05&26.630&2.663&-.034&$.0^3825$&(-.0524,\,.935)&$.0^2164$&\\
&.01&37.109&3.711&-.039&$.0^3588$&(-.0525,\,.953)&$.0^2119$&\\
 \hline
 40&.20&34.844&1.742&-.011&$.0^3338$&(-.0256,\,.947)&$.0^2105$&-.0256\\
&.10&45.995&2.300&-.015&$.0^3255$&(-.0256,\,.959)&$.0^3514$&\\
&.05&57.168&2.858&-.017&$.0^3205$&(-.0256,\,.967)&$.0^3414$&\\
&.01&90.567&4.528&-.020&$.0^3129$&(-.0256,\,.979)&$.0^3263$&\\
 \hline
 $\infty$&.20&$\infty$&1.803&0&0&[0,\,1)&0&0\\
&.10&$\infty$&2.429&0&0&[0,\,1)&0&\\
&.05&$\infty$&3.074&0&0&[0,\,1)&0&\\
&.01&$\infty$&4.628&0&0&[0,\,1)&0&\\
 \hline
\end{tabular}
\end{center}
\caption{Suppose that $p\ge3$, $N=p+2$, and $\del^2$ is sufficiently small, that is, $\del^2\ts<\u m_{p,\alp}\Lam_{2,p,p+2;\alp}^{**}$ ($\del^2\ts<\u m_{p,\alp}'\Lam_{2,p,p+2;\alp}^{**}$). In Case 1 (Case 2) if $\u\psi_{p;\alp}^-<\rho<1$ ($\u\rho_{p;\alp}^-<\rho<\u\rho_{p;\alp}^+)$ then the power of at least one of the ${p\choose2}$ bivariate size-$\alp$ $T^2$-tests dominates that of the overall $p$-variate  size-$\alp$ $T^2$-test. Note that the feasible range of $\rho$ is $\big(-\frac{1}{p-1},1\big)$.}
\label{table:4}
\end{table}
\vskip6pt

\newpage

\nid{\bf 4. Some exact power comparisons for the bivariate case.} The results in Sections 2 and 3 compare the power of the overall $T^2$-test with those of univariate or bivariate $T^2$-tests based on the original variates. However these power comparisons are asymptotic or local, relevant only for noncentrality parameters $\Lam$ that approach $\infty$ or 0. In this section we consider the bivariate case $p=2$ and attempt to compare the exact power functions of the $T^2$-test and the two univariate $t^2$-tests for all values of $\Lam$. Two conjectures are presented; the first of these is confirmed in Proposition 4.3 and applied in Example 4.4 for only two simple cases.
\vskip4pt

\nid{\bf Conjecture 4.1 (weak).} Suppose that $p=2$ and $N$ is odd: $N=2l+1$. Then
for each $\lam>0$, $\exists\, 0<\alp_l^*(\lam)<1$ such that 
\begin{equation}\label{ineq111}
\ts0<\alp<\alp_l^*(\lam)\ \ \implies\ \ \pi_\alp(\lam;\,1,\,2l)>\pi_\alp\big(\big(\frac{4l}{2l-1}\big)\lam;\,2,\,2l-1\big),
\end{equation}
with equality when $\alp=\alp_l^*(\lam)$.\hfill$\square$
\vskip4pt

Conjecture 4.1 is established below for $l=1,2$ and we expect it to hold for all $l\ge3$ as well. However, it is unsatisfactory in that if $\alp_l^*(\lam)$ depends nontrivially on $\lam$ then we cannot conclude that, at least for small $\alp$, one or both of the two univariate size-$\alp$ $t^2$-tests dominate the bivariate size-$\alp$ $T^2$ test in a large region of the alternative hypothesis. For this the following stronger result would be needed.
\vskip4pt

\nid{\bf Conjecture 4.2 (strong).}  Conjecture 4.1 holds with $\alp_l^*(\lam)$ not depending on $\lam$, i.e., $\alp_l^*(\lam)=\alp_l^*$.\hfill$\square$
\vskip4pt

\nid At this time we do not have evidence either for or against Conjecture 4.2. If valid, it would be essential to determine or approximate the values of $\alp_l^*$.
\vskip4pt

%\newpage

\nid{\bf Proposition 4.3.} Conjecture 4.1 is valid for $l=1$ and 2.
\vskip4pt

\nid{\bf Proof.} By \eqref{6},
\begin{equation}\label{piineq1}
\pi_\alp(\lam;\,1,\,2l)>\ts\pi_\alp(2(1+\del)\lam;\,2,\,2l-1)
\end{equation}
if and only if
\begin{align}
s_{\del,\lam}^{(l)}(\alp):=\ts e^{(\frac{1}{2}+\del)\lam}\sum_{k=0}^\infty\ts(\frac{\lam}{2})^k\frac{1}{k!}c_{1,2l;k;\alp}>\sum_{k=0}^\infty\frac{(1+\del)^k\lam^k}{k!}c_{2,2l-1;k;\alp}=:t_{\del,\lam}^{(l)}(\alp).\label{ineq2}
\end{align}
From \eqref{cmnk2} and \eqref{betapdf} we find that
\begin{align}
c_{1,2l;k;\alp}&=\Pr[b_{2l,1+2k}<b_{2l,1;\alp}]\label{newa}\\
&=\ts\frac{\Gam(l+\frac{1}{2}+k)}{\Gam(l)\Gam(\frac{1}{2}+k)}
\int_0^{b_{2l,1;\alp}}b^{l-1}(1-b)^{k-\frac{1}{2}}db;\label{newk}\\ 
\alp&=\Pr[b_{2l,1}<b_{2l,1;\alp}]\label{newb}\\
&=\ts\frac{\Gam(l+\frac{1}{2})}{\Gam(l)\Gam(\frac{1}{2})}
\int_0^{b_{2l,1;\alp}}b^{l-1}(1-b)^{-\frac{1}{2}}db.\label{newc}
\end{align}
Set  $u=1-b$ in \eqref{newb}-\eqref{newc}, then differentiate w.r.to $\alp$ to obtain
\begin{align}
 \ts\frac{d}{d\alp}b_{2l,1;\alp}&=\ts\frac{\Gam(l)\Gam(\frac{1}{2})}{\Gam(l+\frac{1}{2})}\frac {(1-b_{2l,1;\alp})^{\frac{1}{2}} }{b_{2l,1;\alp}^{l-1}};
\label{newee}\\
%\hskip50pt
 \ts\frac{d}{d\alp}c_{1,2l;k,\alp}&=\ts\frac{\Gam(l+\frac{1}{2}+k)\Gam(\f12)}{\Gam(l+\f12)\Gam(\frac{1}{2}+k)}(1-b_{2l,1;\alp})^k\\
%\hskip50pt 
\ts\frac{d}{d\alp}s_{\del,\lam}^{(l)}(\alp)&=\ts e^{(\frac{1}{2}+\del)\lam}\sum_{k=0}^\infty\ts(\frac{\lam}{2})^k\frac{1}{k!}\frac{\Gam(l+\frac{1}{2}+k)\Gam(\f12)}{\Gam(l+\f12)\Gam(\frac{1}{2}+k)}(1-b_{2l,1;\alp})^k\\
&=\ts\big[\sum_{k=0}^\infty\frac{(\frac{1}{2}+\del)^k\lam^k}{k!}\big]\big[\sum_{k=0}^\infty\ts(\frac{\lam}{2})^k\frac{1}{k!}\frac{\Gam(l+\frac{1}{2}+k)\Gam(\f12)}{\Gam(l+\f12)\Gam(\frac{1}{2}+k)}(1-b_{2l,1;\alp})^k\big]\\
&=\ts\sum_{k=0}^\infty\frac{(\frac{1}{2}+\del)^k}{k!}\big[\sum_{r=0}^k{k\choose r}\frac{1}{(1+2\del)^r}\frac{\Gam(l+\frac{1}{2}+r)\Gam(\f12)}{\Gam(l+\f12)\Gam(\frac{1}{2}+r)}(1-b_{2l,1;\alp})^r\big]\lam^k.
\end{align}

Next,
\begin{align}
\alp&=\Pr[b_{2l-1,2}<b_{2l-1,2;\alp}]\label{newg}\\
&=\ts\frac{\Gam(l+\frac{1}{2})}{\Gam(l-\frac{1}{2})\Gam(1)}
\int_0^{b_{2l-1,2;\alp}}b^{l-\frac{3}{2}}db=b_{2l-1,2;\alp}^{l-\f12},\label{newh}\\
b_{2l-1,2;\alp}&=\alp^{\frac{2}{2l-1}};\label{newi}\\
c_{2,2l-1;k;\alp}&=\Pr[b_{2l-1,2+2k}<\alp^{\frac{2}{2l-1}}]\label{newf}\\
&=\ts\frac{\Gam(l+\frac{1}{2}+k)}{\Gam(l-\frac{1}{2})\Gam(1+k)}
\int_0^{\alp^{\frac{2l}{2l-1}}}b^{l-\frac{3}{2}}(1-b)^kdb\label{newj}\\
&=\ts\frac{\Gam(l+\frac{1}{2}+k)} {\Gam(l+\frac{1}{2})\Gam(1+k)}
\int_0^\alp(1-w^{\frac{2l}{2l-1}})^kdw;\label{newkA}\\
\ts\frac{d}{d\alp}c_{2,2l-1;k,\alp}&=\ts\frac{\Gam(l+\frac{1}{2}+k)} {\Gam(l+\frac{1}{2})\Gam(1+k)}
(1-\alp^{\frac{2l}{2l-1}})^k;\label{newl}\\
\ts\frac{d}{d\alp}t_{\del,\lam}^{(l)}(\alp)
&=\ts\sum_{k=0}^\infty\frac{(1+\del)^k\lam^k}{k!}\frac{\Gam(l+\frac{1}{2}+k)} {\Gam(l+\frac{1}{2})\Gam(1+k)}(1-\alp^{\frac{2l}{2l-1}})^k.
\end{align}
Therefore a sufficient condition that $\frac{d}{d\alp}s_{\del,\lam}^{(l)}(\alp)> \frac{d}{d\alp}t_{\del,\lam}^{(l)}(\alp)$ 
%for all $\lam\ge0$ 
is that for all $k\ge0$,
\begin{align}
\ts\frac{(\frac{1}{2}+\del)^k}{(1+\del)^k(1-\alp^{\frac{2l}{2l-1}})^k}\sum_{r=0}^k{k\choose r}\frac{(1-b_{2l,1;\alp})^r}{(1+2\del)^r}\frac{\Gam(l+\frac{1}{2}+r)}{\Gam(\frac{1}{2}+r)}
&\ts\ge\frac{\Gam(l+\frac{1}{2}+k)}{\Gam(\frac{1}{2})\Gam(1+k)},\label{ineqBalp}
\end{align}
with strict inequality for at least one $k$.

Thus for $\alp=0$, a sufficient condition that $\frac{d}{d\alp}s_{\del,\lam}^{(l)}(\alp=0)> \frac{d}{d\alp}t_{\del,\lam}^{(l)}(\alp=0)$ 
%for all $\lam\ge0$ 
is that for all $k\ge0$,
\begin{align}
\ts\frac{(\frac{1}{2}+\del)^k}{(1+\del)^k}\sum_{r=0}^k{k\choose r}\frac{1}{(1+2\del)^r}\frac{\Gam(l+\frac{1}{2}+r)}{\Gam(\frac{1}{2}+r)}
&\ts\ge\frac{\Gam(l+\frac{1}{2}+k)}{\Gam(\frac{1}{2})\Gam(1+k)},
\label{ineqB}
\end{align}
with strict inequality for at least one $k$. After some algebra, \eqref{ineqB} can be written equivalently as
\begin{align}
\ts\E\Big[\frac{\Gam(l+\f12+R_{k,\del})} {\Gam(\f12+R_{k,\del})}\Big]\ge\frac{\Gam(l+\frac{1}{2}+k)}{\Gam(\frac{1}{2})\Gam(1+k)},\label{Binkdel}
\end{align}
where $R_{k,\del}\sim{\rm Binomial}(k,\frac{1}{1+\del})$. Because $R_{k,\del}$ is (strictly) stochastically decreasing in $\del$ (for $k\ge1$) while $\frac{\Gam(l+\f12+R_{k,\del})} {\Gam(\f12+R_{k,\del})}$ is (strictly) increasing in $R_{k,\del}$, the left side of \eqref{Binkdel} $\equiv$ \eqref{ineqB} is (strictly) decreasing in $\del$ (for $k\ge1$).

For $k=0$ both sides of \eqref{ineqB} $=\frac{\Gam(l+\frac{1}{2})}{\Gam(\frac{1}{2})}$. For $k=1$, \eqref{ineqB} is equivalent to the inequality
\begin{align}
\ts\frac{(\frac{1}{2}+\del)}{(1+\del)}\Big[\frac{\Gam(l+\frac{1}{2})}{\Gam(\frac{1}{2})}+\frac{1}{(1+2\del)}\frac{\Gam(l+\frac{3}{2})}{\Gam(\frac{3}{2})}\Big]
\ge&\ \ts\frac{\Gam(l+\frac{3}{2})}{\Gam(\frac{1}{2})},\nonumber
%\ts\del\le&\ \ts\frac{1}{2l-1}.\label{ineqBB}
\end{align}
which is equivalent to $\del\le\frac{1}{2l-1}$. Therefore the sufficient condition \eqref{ineqB} for 
\begin{equation*}
\ts\frac{d}{d\alp}s_{\del,\lam}^{(l)}(\alp=0)> \frac{d}{d\alp}t_{\del,\lam}^{(l)}(\alp=0)
\end{equation*}
will be satisfied for all  $\del\le\frac{1}{2l-1}$ if \eqref{ineqB} holds for $\del=\frac{1}{2l-1}$ for all $k\ge2$, with strict inequality for at least one $k\ge2$.

Because $s_{\del,\lam}^{(l)}(0)=t_{\del,\lam}^{(l)}(0)=0$, it follows from \eqref{piineq1}-\eqref{ineq2} that \eqref{ineqB}, with strict inequality for some $k\ge2$, is a sufficient condition that for each $\lam>0$, $\exists\, 0<\alp_l^*(\lam)<1$ such that 
\begin{equation}\label{ineq11}
\ts0<\alp<\alp_l^*(\lam)\ \ \implies\ \ \pi_\alp(\lam;\,1,\,2l)>\pi_\alp\big(\big(\frac{4l}{2l-1}\big)\lam;\,2,\,2l-1\big),
\end{equation}
with equality when $\alp=\alp_l^*(\lam)$.  

For the simplest case $l=1$ ($N=3$), \eqref{ineqB} with $\del=\frac{1}{2l-1}=1$ becomes
\begin{align}
\ts\big(\frac{3}{4}\big)^k\sum_{r=0}^k{k\choose r}\frac{1}{3^r}(\frac{1}{2}+r)
&\ts\ge\frac{\Gam(\frac{3}{2}+k)}{\Gam(\frac{1}{2})\Gam(1+k)},
\label{ineqBB}
\end{align}
which by \eqref{Binkdel} can be reduced to the equivalent form
\begin{align}
%\ts(\frac{1}{2}+\frac{k}{4})
\ts1+\frac{k}{2}
&\ts\ge\frac{\Gam(\frac{3}{2}+k)}{\Gam(\frac{3}{2})\Gam(1+k)}.
\label{ineqBBB}
\end{align}
It is straightforward to verify \eqref{ineqBBB} by induction on $k$, with strict inequality holding for large $k$ because $\frac{\Gam(\frac{3}{2}+k)}{\Gam(1+k)}=O(k^{\f12})$. Therefore we conclude that 
\eqref{ineq11} holds for $l=1$; that is, for each $\lam>0$, $\exists\, 0<\alp_1^*(\lam)<1$ such that 
\begin{equation}\label{ineq1}
0<\alp<\alp_1^*(\lam)\ \ \implies\ \ \pi_\alp(\lam;\,1,\,2)>\pi_\alp(4\lam;\,2,\,1),
\end{equation}
with equality when $\alp=\alp_1^*(\lam)$. 
\vskip4pt

Next consider the case $l=2$ ($N=5$). With $\del=\frac{1}{2l-1}=\frac{1}{3}$,  \eqref{ineqB}  becomes
\begin{align}
\ts\big(\frac{5}{8}\big)^k\sum_{r=0}^k{k\choose r}\big(\frac{3}{5}\big)^r(\frac{3}{2}+r)(\frac{1}{2}+r)
&\ts\ge\frac{\Gam(\frac{5}{2}+k)}{\Gam(\frac{1}{2})\Gam(1+k)},
\label{ineqCC}
\end{align}
which can be reduced to the equivalent form
\begin{align}
\ts\frac{3}{4}+\frac{9k}{8}+\frac{9k(k-1)}{64}
&\ts\ge \frac{\Gam(\frac{5}{2}+k)}{\Gam(\frac{1}{2})\Gam(1+k)}=(\frac{3}{2}+k)(\frac{1}{2}+k)\frac{\Gam(\frac{1}{2}+k)}{\Gam(\frac{1}{2})\Gam(1+k)}.\label{CCC}
\end{align}
Interestingly, \eqref{CCC} holds with equality for $k=2$ as well as for 0 and 1. Rewrite \eqref{CCC} in the equivalent form
\begin{align}
\ts\frac{k(k-1)\cdots2\cdot1}{(k-\frac{1}{2})(k-\frac{1}{2})\cdots\frac{3}{2}\cdot\frac{1}{2}}
&\ts\ge\frac{48+128k+64k^2}{48+63k+9k^2}.\label{ineqN}
\end{align}
To verify \eqref{ineqN} by induction on $k$, it suffices to show that for $k\ge2$,
\begin{align}
\ts\frac{k+1}{k+\f12}\frac{48+128k+64k^2}{48+63k+9k^2}\ge\frac{48+128(k+1)+64(k+1)^2}{48+63(k+1)+9(k+1)^2}.\label{ineqP}
\end{align}
After simplification, this is equivalent to the inequality
\begin{align}
4k^3+4k^2-5k-3\ge0,\label{ineqQ}
\end{align}
which holds for all $k\ge1$, with strict inequality for $k+1\ge3$. %Therefore 
%$\frac{d}{d\alp}s_{\del,\lam}^{(2)}(\alp=0)> \frac{d}{d\alp}t_{\del,\lam}^{(2)}(\alp=0)$ holds for  $\del\le\frac{1}{3}$. This implies that 
Therefore we conclude that \eqref{ineq11} holds for $l=2$:
for each $\lam>0$, $\exists\, 0<\alp_2^*(\lam)<1$ such that 
\begin{equation}\label{ineq177}
\ts0<\alp<\alp_2^*(\lam)\ \ \implies\ \ \pi_\alp(\lam;\,1,\,4)>\pi_\alp(\frac{8}{3}\lam;\,2,\,3),
\end{equation}
with equality when $\alp=\alp_2^*(\lam)$. \hfill$\square$
\vskip6pt

\nid{\bf Example 4.4.} Return to the bivariate Example 2.1, where $p=2$ and 
\begin{equation}\label{lams5}
\ts\Lam_1=N\gam_1^2,\quad\Lam_2=N\eta^2\gam_1^2,\quad \Lam=N\big(\frac{1-2\eta\rho+\eta^2}{1-\rho^2}\big)\gam_1^2;
\end{equation}
(recall \eqref{lams}). For $N=3$ it follows from \eqref{ineq1} and \eqref{lams5} that for each $\gam_1^2>0$,
\begin{equation}\label{ineq15}
0<\alp<\alp_1^*(3\gam_1^2)\ \ \implies\ \ \pi_\alp(3\gam_1^2;\,1,\,2)>\pi_\alp(4\!\cdot\!3\gam_1^2;\,2,\,1).
\end{equation}
%with equality when $\alp=\alp_1^*(\lam)$. 
Furthermore,
\begin{align}
\ts4\!\cdot\!3\gam_1^2>3\big(\frac{1-2\eta\rho+\eta^2}{1-\rho^2}\big)\gam_1^2
\ \ &\iff\ \ \ts4>\frac{1-2\eta\rho+\eta^2}{1-\rho^2}\label{ineq16}\\
&\iff0>4\rho^2-2\eta\rho+(\eta^2-3):=h_\eta(\rho).\label{ineq17}
\end{align}
The two roots of $h_\eta(\rho)$ are $\ts\^\rho_\eta^\pm=\frac{\eta\pm\sqrt{12-3\eta^2}}{4}$; note that $\^\rho_{-\eta}^{\pm}=-\^\rho_{\eta}^{\mp}$. Some values appear in Table 5. Thus, if $\alp<\alp_1^*(3\gam_1^2)$ then
\begin{align}
\rho\in(\^\rho_{\eta}^-,\,\^\rho_{\eta}^+)&\implies\max(\pi_\alp(\Lam_1;\,1,2),\,\pi_\alp(\Lam_2;\,1,2))>\pi_\alp(\Lam;\,2,1)\label{rhoetapm}\\
&\implies |\^\ome_\alp(\gam_1,\eta,\rho)|=1;\label{rhoetapm1}
\end{align}
that is, at least one of the two univariate $t^2$-tests is more powerful than the overall bivariate $T^2$-test. This occurs in the $(\eta,\rho)$-regions of the parameter space indicated in Table 5, which constitute a substantial part of the alternative hypothesis. 

Similarly, for $N=5$ it follows from \eqref{ineq177} and \eqref{lams5} that for each $\gam_1^2>0$,
\begin{equation}\label{ineq155}
0<\alp<\alp_2^*(5\gam_1^2)\ \ \implies\ \ \pi_\alp(5\gam_1^2;\,1,\,4)>\ts\pi_\alp(\frac{8}{3}\!\cdot\!5\gam_1^2;\,2,\,3).
\end{equation}
%with equality when $\alp=\alp_1^*(\lam)$. 
Furthermore,
\begin{align}
\ts\frac{8}{3}\!\cdot\!5\gam_1^2>5\big(\frac{1-2\eta\rho+\eta^2}{1-\rho^2}\big)\gam_1^2
\ \ &\iff\ \ \ts\frac{8}{3}>\frac{1-2\eta\rho+\eta^2}{1-\rho^2}\label{ineq161}\\
&\iff0>\ts\frac{8}{3}\rho^2-2\eta\rho+(\eta^2-\frac{5}{3}):=\~h_\eta(\rho).\label{ineq178}
\end{align}
The two roots of $\~h_\eta(\rho)$ are $\ts\~\rho_\eta^\pm=\frac{3\eta\pm\sqrt{40-15\eta^2}}{8}$; note that $\~\rho_{-\eta}^{\pm}=-\~\rho_{\eta}^{\mp}$. Some values appear in Table 5. Thus, if $\alp<\alp_2^*(5\gam_1^2)$ then
\begin{align}
\rho\in(\~\rho_{\eta}^-,\,\~\rho_{\eta}^+)&\implies\max(\pi_\alp(\Lam_1;\,1,4),\,\pi_\alp(\Lam_2;\,1,4))>\pi_\alp(\Lam;\,2,3)\label{rhoetapm2}\\
&\implies |\^\ome_\alp(\gam_1,\eta,\rho)|=1;\label{rhoetapm3}
\end{align}
that is, at least one of the two univariate $t^2$-tests is more powerful than the overall bivariate $T^2$-test. Again this occurs in the $(\eta,\rho)$-regions of the parameter space indicated in Table 5. 

Thus for $p=2$, $N=3$ or 5, and sufficiently small $\alp$ (but depending on $\gam_1^2$), the overall size-$\alp$ $T^2$-test is dominated by at least one of the two univariate size-$\alp$ $t^2$-tests over a fairly large portion of the alternative hypothesis.\hfill$\square$

\begin{table}[h!]
\centering
\begin{center}
\begin{tabular}{ |c|c|c|c|c|c|c|c|c| }
%\begin{tabular}{ |p{2.4cm}||p{.6cm}|p{.6cm}|p{.6cm}|p{.6cm}|p{.6cm}|p{.6cm}|p{.6cm}|p{.6cm}|p{.6cm}|p{.6cm}|p{.6cm}| }
% \hline
%q&\|a^{(qœ§ )}\|&\^y^{(q)}'&y^{(q)}'&3&4&5&6&7&8&9&10\\
 \hline
 $\eta$&$(\^\rho_{\eta}^-,\,\^\rho_{\eta}^+)$&$(\~\rho_{\eta}^-,\,\~\rho_{\eta}^+)$\\
 \hline
 \hline
 1&(-.5,\,1)&(-.25,\,1)\\
 .75&(-.615,\,.990)&(-.421,\,.983)\\
.5&(-.714,\,.964)& (-.566,\,.940)\\
 .25&(-.797,\,.922)&(-.687,\,.875)\\
0&(-.866,\,.866)&(-.791,\,.791)\\
 \hline
\end{tabular}
\end{center}
\caption{For $p=2$, $N=3$ ($N=5$), and $\alp<\alp_1^*(3\gam_1^2)$ ($\alp<\alp_2^*(5\gam_1^2)$),
if $\rho\in(\^\rho_{\eta}^-,\,\^\rho_{\eta}^+)$ ($\rho\in(\~\rho_{\eta}^-,\,\~\rho_{\eta}^+)$) then $ |\^\ome_\alp(\gam_1,\eta,\rho)|=1$, so the power of at least one of the two univariate size-$\alp$ $t^2$-tests exceeds that of the bivariate size-$\alp$  $T^2$-test. (Note that $\^\rho_{-\eta}^{\pm}=-\^\rho_{\eta}^{\mp}$ and $\~\rho_{-\eta}^{\pm}=-\~\rho_{\eta}^{\mp}$.)}
\label{table:5}
\end{table}

%COMMENT
\begin{comment}
\nid{\bf Example 4.4.NOT USED} Return to Example 2.2, where $N\ge5$ and $p=N-2\ge3$. From \eqref{ineq1},
\begin{equation}\label{alpalppipi}
\ts0<\alp<\alp_2^*(\Lam^{(2)})\ \ \implies\ \ \pi_\alp(\Lam^{(2)};\,1,\,4)>\pi_\alp(\frac{8}{3}\Lam^{(2)};\,2,\,3).
\end{equation}
{\bf BAD BECAUSE $\alp_2^*(\Lam^{(2)})$ depends on RHO???}

In Case 1 (recall \eqref{Lam2lam})
\begin{align}
\ts\Lam^{(2)}=\frac{2N\del^2}{1+\rho},&\qquad\ \ \ \ \ \Lam=\ts\frac{(N-2)N\del^2}{1+(N-3)\rho},\label{Lam2lamA}\\
\ts\frac{8}{3}\Lam^{(2)}\ge\Lam\ \ &\iff\ \ \rho>\ts\frac{3N-22}{13N-42}=:\~\psi_N^-.\label{83psitilde}
\end{align}
Thus if $\alp<\alp_2^*(\Lam^{(2)})$ then....

Case 2

Case 3

Case 4
\end{comment}
%END

\vskip6pt

\nid{\bf 5. Concluding remarks.} For the purpose of stimulating future research, the questions raised in this report are stated formally as follows:
\vskip4pt

\nid {\it The Oracular Variable-Selection Problem (OVSP)} is that of determining the function $\^\ome_\alp(\gam,R)$, as defined in \eqref{omestar}, and using this to determine
 the regions
\begin{align}
%&\hskip120pt 
A_\alp(i)\equiv\big\{(\gam,R)\bigm| &|\^\ome_\alp(\gam,R)|=i\big\},\quad i=1,...,p.\label{optreg}
%\\ =\{(\gam,R)\mid \pi_\alp(\Lam_\ome(\gam,R);\,|\ome|,\,N-|\ome|)&\ge\pi_\alp(\Lam_{\ome'}(\gam,R) ;\,|\ome'|,\,N-|\ome'|)\ \forall\ome'\in\Ome_p\}.\nonumber
\end{align}
%(Note that these regions may overlap.) 
\vskip4pt
\nid \nid {\it The Parsimonious Variable-Selection Problem (PVSP)} asks if $A_\alp(i)$ comprises a substantial portion of the alternative hypothesis $K$ for small values of $i$, e.g., $i=1,2$.
\vskip4pt

\nid If the answer to the PVSP is positive, then variable selection in some applied investigations can be limited to small, easily interpretable subsets of variables.
\vskip6pt

\nid{\bf Acknowledgment.} Special thanks are due to David Perlman for raising the questions addressed here and providing some supporting data.
\vskip24pt

%\newpage

\centerline{\bf References}
\vskip6pt

\nid [A] Anderson, T. W. (2003). {\it An Introduction to Multivariate Statistical Analysis, 3rd edition,} Wiley \& Sons, New York.
\smallskip

\nid [DGP] Das Gupta, S. and Perlman, M. D. (1974). On the power of the noncentral $F$-test: effect of additional variates on Hotelling's $T^2$-test.  {\it J.Amer. Statist. Assoc.} {\bf 69} 174-180. 

\smallskip

\nid [GK] Giri, N. and Kiefer, J. (1964).  Local and asymptotic minimax properties of multivariate tests. {\it Ann. Math. Statist.} {\bf 35} 21-35.

\smallskip

\nid [GKS] Giri, N., Kiefer, J., and Stein, C. (1963). Minimax character of Hotelling's $T^2$-test in the simplest case. {\it Ann. Math. Statist.} {\bf 34} 1524-1535.

\smallskip

\nid [KS] Kiefer, J. and Schwartz, R. (1965). Admissible Bayes character of $T^2$-, $R^2$, and other fully invariant tests for classical multivariate normal testing problems. {\it Ann. Math. Statist.} {\bf 36} 747-770.

\smallskip

%\newpage

\nid [MP] Marden, J. and Perlman, M. D. (1980). Invariant tests for means with covariates. {\it Ann. Statist.} {\bf 8} 25-63.

\smallskip

\nid [S] Stein, C. (1956). The admissibility of Hotelling's $T^2$-test.  {\it Ann. Math. Statist.} {\bf 27} 616-623.
\vskip24pt

\centerline{\bf Appendix}
\vskip4pt

\nid{\bf A. Testing for additional information.} Variable selection for the $T^2$-test and related linear discriminant analysis was thoroughly studied in the 1970s and 1980s, an era of limited computer power, and subsequently by several authors with greater ability to consider all-subsets methods; a list of references appears below. Almost all of these studies were based on testing for additional information (= increased Mahalanobis distance),
%Rao's test or Subrahmaniam's test for additional information, 
as now described.

For any two nested subsets $\ome\subset\ome'$ in $\Ome_p$, in general $\Lam_\ome\le\Lam_{\ome'}$. The question of whether the power of the $T_{\ome'}$-test exceeds that of the $T_{\ome}$-test for the testing problem \eqref{H0K} usually was formulated as the problem of {\it testing for additional information (TAI),} namely, testing
\begin{equation}\label{TAI}
\Lam_\ome'=\Lam_{\ome}\quad \mathrm{vs.}\quad \Lam_\ome'>\Lam_{\ome}
\end{equation}
% the likelihood ratio test (LRT) for testing $\Lam_{\ome'}=\Lam_\ome$ vs. $\Lam_{\ome'}>\Lam_\ome$, 
based on a preliminary sample -- see [R] \S8c.4. This formulation was adopted by many researchers, even while citing the following result of [DGP] which implies that this standard formulation of TAI is inappropriate.

It was shown in [DGP, Theorem 2.1] that for fixed $\lam>0$, the power function $\pi_\alp(\lam;m,n)$ (recall \eqref{pialp}) of the non-central $f$-test is strictly decreasing in $m$ and strictly increasing in $n$.\footnote{\label{Foot1} 
%%In the display following (3.10) on p.177 of [DGP], the upper limit of summation should be $b_i-1$. 
In line 2 of the second column on p.179 of [DGP], ``conclude" should be``include". In the line following the third display in the second column on p.179, ``$j$" should be``$f$". In Remark 4.1 on p.180, ``increasing in $m$" should be ``decreasing in $m$". In the next line, ``$m\to\infty$" should be ``$n\to\infty$".} Therefore for any integer $1\le q\le n-1$ there exists a unique real number
\begin{equation}\label{7}
g_\alp(\lam)\equiv g_\alp(\lam; m,n,q)>0
\end{equation}
such that
\begin{equation}\label{8}
\pi_\alp(\lam;m,n)=\pi_\alp(\lam+g_\alp(\lam);m+q, n-q).
\end{equation}
Here $g_\alp(0)=0$ and $g_\alp(\lam)$ is strictly increasing in $\lam$; cf. [DGP, Theorem 3.1]. Thus the power is increased only if
\begin{equation}\label{lamgrlam}
\Lam_{\ome'}>\Lam_\ome+g_\alp(\Lam_{\ome}; |\ome|,N-|\ome|,|\ome'\stms\ome|).
\end{equation}
Therefore [DGP, Section 4] introduced the problem of {\it testing for increased power (TIP),} namely, testing
\begin{align*}
&H_1:\Lam_{\ome'}\le\Lam_\ome+g_\alp(\Lam_{|\ome|}; |\ome|,N-|\ome|,|\ome'\stms\ome|)\\
 \mathrm{vs.}\ \ & K_1:\Lam_{\ome'}>\Lam_\ome+g_\alp(\Lam_{|\ome|}; |\ome|,N-|\ome|,|\ome'\stms\ome|),
\end{align*}
and proposed several (approximate) tests. 

This proposal was noted by subsequent authors but never implemented for variable selection, possibly because of difficulties in computing the functions $g_\alp(\cdot)$, especially if many pairs $(\ome,\ome')$ must be considered. However, if as suggested above, variable selection might be limited to very small subsets of variables in practical applications, then replacing the TAI by the TIP might be feasible.
\vskip4pt

\nid{\bf Remark A.1.} The relation \eqref{ineq111} in Conjecture 4.1 can be stated equivalently in terms of $g_\alp$:
\begin{equation}\label{g12ineq1111}
\ts 0<\alp<\alp_l^*(\lam)\ \ \implies\ \ g_\alp(\lam;1,2l,1)>\big(\frac{2l+1}{2l-1}\big)\lam\quad\forall\ \lam>0,
\end{equation}
with equality when $\alp=\alp_l^*(\lam)$. Thus the relations \eqref{ineq1} and \eqref{ineq177} in Proposition 4.3 also can be stated equivalently in terms of $g_\alp$:
\begin{equation}\label{g12ineqA}
0<\alp<\alp_1^*(\lam)\ \ \implies\ \ g_\alp(\lam;1,2,1)>3\lam\quad\forall\ \lam>0,
\end{equation}
with equality when $\alp=\alp_1^*(\lam)$;
\begin{equation}\label{g12ineqAB}
\ts 0<\alp<\alp_2^*(\lam)\ \ \implies\ \ g_\alp(\lam;1,4,1)>\frac{5}{3}\lam\quad\forall\ \lam>0,
\end{equation}
with equality when $\alp=\alp_2^*(\lam)$.
\vskip24pt

\centerline{\bf Additional References for the Appendix}
\vskip6pt

\nid [H1] Hand, D. J. (1981). {\it Discrimination and Classification,} Wiley \& Sons, New York. [Chapter 6]
\smallskip

\nid [H2] Hawkins, D. M. (1976). The subset problem in multivariate analysis of variance. {\it J. Royal Statist. Soc, Series B} {\bf 38} 132-139.
\smallskip

\nid [JW] Jain, A. K. and Waller, W. G. (1978). On the optimal number of features in the classification of multivariate Gaussian data. {\it Pattern Recognition} {\bf 10} 103-109.
\smallskip

\nid [JWT] Jiang, W., Wang, K., and Tsung, F. (2012). A variable-selection-based multivariate EWMA chart for process monitoring and diagnosis. {\it J. Quality Technology} {\bf 44} 209-230.
\smallskip

\nid [McC] McCabe, G. P. Jr. (1975). Computations for variable selection in discriminant analysis. {\it Technometrics} {\bf 17} 259-263.
\smallskip

\nid [McK] McKay, R. J. (1976). A graphical aid to selection of variables in two-group discriminant analysis. {\it Appl. Statist.} {\bf 27} 259-263.
\smallskip

\nid [McL1] McLachlan, G. J. (1976). On the relationship between the $F$ test and the overall error rate for variable selection in two-group discriminant function. {\it Biometrics} {\bf 36} 501-510.
\smallskip

\nid [McL2] McLachlan, G. J. (1980). A criterion for selecting variables for the linear discriminant function. {\it Biometrics} {\bf 32} 529-534.
\smallskip

\nid [McL3] McLachlan, G. J. (1992). {\it Discriminant Analysis and Statistical Pattern Recognition,} Wiley \& Sons, New York. [Chapter 12]
\smallskip

\nid [Mu] Murray, G. D. (1977). A cautionary note on selection of variables in discriminant analysis. {\it Appl. Statist.} {\bf 26} 246-250.
\smallskip

\nid [NT] Nobuo, S. and Takahisa, I. (2016). A variable selection method for detecting abnormality base on the $T^2$ test. {\it Comm. Statist. - Theory and Methods} {\bf 46} 501-510.
\smallskip

\nid [R] Rao, C. R. (1973). {\it Linear Statistical Inference and its Applications, 2nd edition,} Wiley \& Sons, New York.
\smallskip

\nid [S] Schaafsma, W. (1982). In {\it Handbook of Statistics, Vol. 2,} P. R. Krishnaiah and L. N. Kanal, eds., 857-881. North Holland, Amsterdam.
\smallskip

\end{document}